\documentclass[leqno,12pt,twoside]{amsart}
\usepackage{geometry}

\usepackage{times}

\usepackage{amsmath, amssymb, amsfonts, latexsym, mdwlist, amsthm}
\usepackage{subfig}
\usepackage{graphicx}
\usepackage{wrapfig}

\usepackage{tikz}
\usetikzlibrary{calc,trees,positioning,arrows,chains,shapes.geometric,%
    decorations.pathreplacing,decorations.pathmorphing,shapes,%
    matrix,shapes.symbols}

\tikzset{
>=stealth',
  punktchain/.style={
    rectangle,
    rounded corners,
    draw=black, thick,
    minimum height=3em,
    text centered,
    on chain},
  line/.style={draw, thick, <-},
  element/.style={
    tape,
    top color=white,
    bottom color=blue!50!black!60!,
    minimum width=8em,
    draw=blue!40!black!90, very thick,
    text width=10em,
    minimum height=3.5em,
    text centered,
    on chain},
  every join/.style={->, thick,shorten >=1pt},
  decoration={brace},
  tuborg/.style={decorate},
  tubnode/.style={midway, right=2pt},
}

\usepackage{paralist}
\setdefaultenum{(a)}{(i)}{}{}
\usepackage{enumitem} 

\def\C{\ensuremath{\mathbb{C}}}
\def\D{\ensuremath{\mathbb{D}}}

\def\N{\ensuremath{\mathbb{N}}}
\def\P{\ensuremath{\mathbb{P}}}
\def\Q{\ensuremath{\mathbb{Q}}}
\def\R{\ensuremath{\mathbb{R}}}
\def\Z{\ensuremath{\mathbb{Z}}}

\def\Aut{\mathop{\mathrm{Aut}}\nolimits}

\def\ch{\mathop{\mathrm{ch}}\nolimits}

\def\Coh{\mathop{\mathrm{Coh}}\nolimits}

\def\cok{\mathop{\mathrm{cok}}}

\def\dim{\mathop{\mathrm{dim}}\nolimits}

\def\inf{\mathop{\mathrm{inf}}\nolimits}

\def\GL{\mathop{\mathrm{GL}}}
\def\Hom{\mathop{\mathrm{Hom}}\nolimits}

\def\RlHom{\mathop{\mathbf{R}\mathcal Hom}\nolimits}

\def\im{\mathop{\mathrm{im}}\nolimits}

\def\min{\mathop{\mathrm{min}}\nolimits}

\def\Num{\mathop{\mathrm{Num}}\nolimits}
\def\NS{\mathop{\mathrm{NS}}\nolimits}

\def\rk{\mathop{\mathrm{rk}}}

\def\td{\mathop{\mathrm{td}}\nolimits}

\newenvironment{Prf}{\textit{Proof.}\/}{\hfill$\Box$}

\def\MG13{\ensuremath{{\mathcal M}_{\Gamma_1(3)}}}
\def\tildeMG13{\ensuremath{\widetilde{\mathcal M}_{\Gamma_1(3)}}}
\def\Stab{\mathop{\mathrm{Stab}}}
\def\into{\ensuremath{\hookrightarrow}}
\def\onto{\ensuremath{\twoheadrightarrow}}

\def\blank{\underline{\hphantom{A}}}

\def\Db{\mathrm{D}^{\mathrm{b}}}

\def\abs#1{\left\lvert#1\right\rvert}

\newcommand\stv[2]{\left\{#1\,\colon\,#2\right\}}

\makeatletter
\newtheorem*{rep@theorem}{\rep@title}
\newcommand{\newreptheorem}[2]{%
\newenvironment{rep#1}[1]{%
 \def\rep@title{#2 \ref{##1}}%
 \begin{rep@theorem}}%
 {\end{rep@theorem}}}
\makeatother

\newtheorem{Thm}{Theorem}[subsection]
\newreptheorem{Thm}{Theorem}
\newtheorem{Thm-s}[Thm]{Theorem}
\newtheorem{Prop-s}[Thm]{Proposition}
\newtheorem{Lem-s}[Thm]{Lemma}
\newtheorem{Prop}[Thm]{Proposition}
\newtheorem{PropDef}[Thm]{Proposition and Definition}
\newtheorem{Lem}[Thm]{Lemma}
\newtheorem{Cor}[Thm]{Corollary}
\newreptheorem{Cor}{Corollary}
\newtheorem{Con}[Thm]{Conjecture}
\newreptheorem{Con}{Conjecture}

\newtheorem{thm-int}{Theorem}

\theoremstyle{definition}
\newtheorem{Def-s}[Thm]{Definition}
\newtheorem{Def}[Thm]{Definition}
\newtheorem{Rem}[Thm]{Remark}

\newtheorem{Prob}[Thm]{Problem}
\newtheorem{Ex}[Thm]{Example}

\def\C{\ensuremath{\mathbb{C}}}
\def\D{\ensuremath{\mathbb{D}}}

\def\N{\ensuremath{\mathbb{N}}}
\def\P{\ensuremath{\mathbb{P}}}
\def\Q{\ensuremath{\mathbb{Q}}}
\def\R{\ensuremath{\mathbb{R}}}
\def\Z{\ensuremath{\mathbb{Z}}}

\def\AA{\ensuremath{\mathcal A}}
\def\BB{\ensuremath{\mathcal B}}
\def\CC{\ensuremath{\mathcal C}}
\def\DD{\ensuremath{\mathcal D}}

\def\FF{\ensuremath{\mathcal F}}

\def\NN{\ensuremath{\mathcal N}}
\def\OO{\ensuremath{\mathcal O}}
\def\PP{\ensuremath{\mathcal P}}

\def\QQ{\ensuremath{\mathcal Q}}
\def\TT{\ensuremath{\mathcal T}}

\def\OPP{\ensuremath{\widehat{\mathcal P}}}

\def\cht{\ensuremath{\ch^B}}

\def\DDD{\mathfrak D}

\newcommand{\ignore}[1]{}

\begin{document}
\title[Stability conditions and Bogomolov-Gieseker type inequalities]{Bridgeland Stability conditions
on threefolds I: Bogomolov-Gieseker type inequalities}

\author{Arend Bayer}
\address{Department of Mathematics, University of Connecticut U-3009, 196 Auditorium Road,
Storrs, CT 06269-3009, USA}
\email{bayer@math.uconn.edu}
\urladdr{http://www.math.uconn.edu/~bayer/}

\author{Emanuele Macr\`i}
\address{Mathematical Institute, University of Bonn, Endenicher Allee 60, D-53115 Bonn, Germany \& Department of Mathematics, University of Utah, 155 S 1400 E, Salt Lake City, UT 84112, USA}
\curraddr{Department of Mathematics, The Ohio State University, 231 W 18th Avenue, Columbus, OH 43210, USA}
\email{macri.6@math.osu.edu}
\urladdr{http://www.math.osu.edu/~macri.6/}

\author{Yukinobu Toda}
\address{Institute for the Physics and Mathematics of the Universe, University of Tokyo, 5-1-5 Kashiwanoha, Kashiwa, 277-8583, Japan}
\email{yukinobu.toda@ipmu.jp}


\date{\today}

\begin{abstract}
We construct new t-structures on the derived category of coherent sheaves on smooth projective threefolds.
We conjecture that they give Bridgeland stability conditions near the large volume limit.
We show that this conjecture is equivalent to a Bogomolov-Gieseker type inequality for the
third Chern character of certain stable complexes.
We also conjecture a stronger inequality, and prove it in the case of projective
space, and for various examples.

Finally, we prove a version of the classical Bogomolov-Gieseker inequality, not involving the
third Chern character, for stable complexes.
\end{abstract}

\maketitle

\setcounter{tocdepth}{1}
\tableofcontents

\section{Introduction}\label{sec:intro}

In this paper, we give a conjectural construction of Bridgeland stability conditions on the derived
category of any projective threefold $X$.  The construction depends on a conjectural
Bogomolov-Gieseker type inequality for objects in the derived category that are stable with respect
to ``tilt-stability'', which is an auxiliary stability condition for two-term complexes on $X$.

\subsection{The existence problem}
Spaces of Bridgeland stability conditions have turned out to be extremely interesting. However, we
do not know a single example of a Bridgeland stability condition on a projective Calabi-Yau
threefold, which is likely to be the most interesting case. The main obstacle is the failure to
solve the following question:

\begin{Prob}
Given a projective threefold $X$, find a heart $\AA \subset \Db(X)$ of a bounded t-structure,
and a group homomorphism  (called central charge)
$Z \colon K(\Db(X)) \to \C$ defined over $\Q$, such that\footnote{The assumption that $Z$ is defined
over $\Q$ is necessary, in practice, to prove the existence of
Harder-Narasimhan filtrations for the induced notion of $Z$-stability on $\AA$.}
\begin{equation} \label{eq:Z-positivity}
0 \neq E \in \AA \quad \Longrightarrow \quad Z(E) \in \stv{re^{i \phi}}{ r > 0, 0 < \phi \le 1}.
\end{equation}
\end{Prob}
We will restrict our attention to central charges $Z$ that are ``numerical'',
i.e., $Z$ factors via the Chern character map $\ch \colon K(\Db(X)) \to \Num_\Q(X)$ to the group 
$\Num_\Q(X)$ of cycles up to numerical equivalence, tensored by $\Q$.

Condition \eqref{eq:Z-positivity} is a highly non-trivial positivity property. For example, it cannot
possible be satisfied when $X$ is projective of dimension $\ge 2$, and $\AA = \Coh X$ is the heart of the
standard t-structure (cf.~\cite[Lemma~2.7]{Toda:limit-stable}).
Further, the construction of
stability conditions for surfaces (see \cite{Bridgeland:K3, Aaron-Daniele}) 
needs the Bogomolov-Gieseker inequality for slope-stable
bundles and the Hodge Index theorem. The methods of \cite{localP2} imply an even closer
relationship: knowing the set of possible numerical central charges
$Z$ for which skyscraper sheaves of points $k(x)$ are stable is essentially \emph{equivalent} to knowing 
the set of Chern characters of slope-semistable bundles for any polarization of $X$.

Motivated by the construction of $\pi$-stability in string theory (see e.g.~
\cite{Aspinwall:Dbranes-CY}), and by the case of curves and surfaces, one can be even
more precise.
Given an ample class $\omega \in \NS_\Q(X)$
and a ``B-field'' $B \in \NS_\Q(X)$, we define a central charge $Z_{\omega, B}$ by
\begin{align} \label{3fold-charge}
Z_{\omega, B}(E) &=-\int_{X}e^{-B-i\omega}\ch(E) \\\nonumber
&= \left(-\cht_3(E)+\frac{\omega^2}{2}\cht_1(E) \right)
+i\left(\omega \cht_2(E)-\frac{\omega^3}{6}\cht_0(E)  \right).
\end{align}
where $\cht$ denotes the twisted Chern character $\cht(E) = e^{-B} \ch(E)$.

\begin{Con}\label{Con:LV}
There exists a heart $\AA_{\omega, B} \subset \Db(X)$ of a bounded t-structure, such that the pair
$(Z_{\omega}, \AA_{\omega})$ is a stability condition on $\Db(X)$ for which skyscraper sheaves
$k(x)$ of points $x \in X$ are stable.
\end{Con}
For $\omega = m \omega_0$ and $m \gg 0$, these would be stability conditions near the ``large-volume
limit''.
As indicated above, the corresponding conjecture is known when $\dim X \le 2$.
In fact when $\dim X=1$, we can take $\AA_{\omega, B}$ to be $\Coh X$. 
When $\dim X=2$, we need to tilt (cf.~Section~\ref{subsec:Tilt}) the abelian category $\Coh X$ to
construct $\AA_{\omega}$. 
We will recall its construction in Proposition~\ref{Rem:BB}.
On the other hand, with very few exceptions (varieties for which $\Db(X)$ admits a complete
exceptional collection), the above conjecture is still open in higher dimension.

\subsection{Our approach}
Given $\omega, B$ as above, we construct a candidate $\AA_{\omega, B}$ for Conjecture
\ref{Con:LV} as a double tilt of $\Coh X$:
\begin{itemize}
\item We first use classical slope-stability
on $\Coh X$ to define $\BB_{\omega, B}$ as a \emph{tilt} of $\Coh X$ with respect to a
torsion pair.
\item We define an analogue of slope-stability on $\BB_{\omega, B}$, and then similarly
define $\AA_{\omega, B}$ as a tilt of $\BB_{\omega, B}$.
\end{itemize}
We now give a sketch of our construction; 
the details will be given in Section \ref{sec:Yukinobu}. Consider the classical slope-stability with
respect to the polarization $\omega$, twisted by $B$: here the slope of a sheaf $\FF$ is given by
$\mu_{\omega, B}(\FF) = \frac{\omega^2 \cht_1(\FF)}{\omega^3\ch_0(\FF)}$. Let
$\TT_{\omega, B} \subset \Coh X$ be the category generated, via extensions, by 
slope-semistable sheaves of slope $\mu_{\omega, B} > 0$ (where torsion sheaves are considered
slope-semistable of slope $\mu_{\omega, B} = +\infty$); and similarly, 
let $\FF_{\omega, B}$ the
subcategory generated by slope-semistable sheaves of slope $\mu_{\omega, B} \le 0$. Following
Bridgeland's construction for K3 surfaces in \cite{Bridgeland:K3}, we define $\BB_{\omega, B} \subset \Db(X)$ by
\[
\BB_{\omega, B} := \langle \FF_{\omega, B}[1], \TT_{\omega, B} \rangle
\]
where $[1]$ denotes the shift and $\langle \blank \rangle$ the extension-closure; see
Section \ref{subsec:Tilt} for alternative descriptions of the tilted heart.

We then define the following slope-function on $\BB_{\omega, B}$:
\[
\nu_{\omega, B}(E)=\frac{\Im Z_{\omega, B}(E)}{\omega^2 \cht_1(E)}
= \frac{\omega \cht_2(E)-\frac{\omega^3}{6}\cht_0(E)} {\omega^2 \cht_1(E)}.
\]
We show that this produces a notion of slope-stability on $\BB_{\omega, B}$, which we call
\emph{tilt-stability}. Using tilt-stability, we can define a torsion pair
$\TT'_{\omega, B}, \FF'_{\omega, B}$ in the category $\BB_{\omega, B}$ exactly as in the
case of slope-stability for $\Coh X$ above.
Tilting at this torsion pair produces a heart $\AA_{\omega, B}$.

We also give a second construction of the same heart in Section \ref{sec:polynomial}, starting
from a category $\Coh^p$ of perverse coherent sheaves, and using polynomial stability conditions 
rather than slope-stability. It is less concrete, but is inherently well-behaved with respect to 
the derived dualizing functor $\RlHom(\blank, \OO_X)$. In Section \ref{sec:Comparison}, we show
that the two constructions agree.

\subsection{Conjectures and Results}
We propose the following conjecture.
\begin{repCon}{Con:stability}
Suppose that $X$ is a smooth projective threefold over $\mathbb{C}$. 
Then the pair $(Z_{\omega, B}, \AA_{\omega, B})$ is a stability condition 
on $\Db(X)$. 
\end{repCon}

At this moment we are not able to show the above conjecture when $X$ is a Calabi-Yau threefold.
As a first evidence for the conjecture, we prove:
\begin{repThm}{thm:P3}
Conjecture~\ref{Con:stability} holds for $X = \P^3$, $B= 0$ and $\omega^3<3\sqrt{3}$.
\end{repThm}
Our method also works for other threefolds with complete exceptional collections.

By construction of $\AA_{\omega, B}$, it is immediate that
$\Im Z_{\omega, B}(E) \ge 0$ for any $E \in \AA_{\omega, B}$. Thus, to show that
condition \eqref{eq:Z-positivity} holds, we only need to consider objects with
$Z_{\omega, B}(E) \in \R$, and have to show that in fact
$Z_{\omega, B}(E) < 0$ in this case.  As in the case of surfaces, this comes down to a
Bogomolov-Gieseker type inequality for tilt-stable objects in $\BB_{\omega, B}$:

\begin{repCon}{Con:stability2}
For any tilt-stable object $E\in \BB_{\omega, B}$ satisfying 
$\nu_{\omega, B}(E) = 0$, i.e., 
\begin{align*}
\frac{\omega^3}{6} \cht_0(E) & =\omega \cht_2(E), \\
\intertext{we have the following inequality:}
\cht_3(E) & < \frac{\omega^2}{2}\cht_1(E). 
\end{align*}
\end{repCon}
In fact, with Corollary \ref{cor:HNfiltrations} we show that Conjecture~\ref{Con:stability} and
Conjecture \ref{Con:stability2} are equivalent.  The essential ingredient is Proposition
\ref{prop:Noether}, which shows that the abelian category $\AA_{\omega, B}$ is Noetherian.

Such a strong Bogomolov-Gieseker type inequality for $\ch_3$ is not available for
slope-semistable sheaves; in fact even for $\P^3$, the best possible results are much worse
(see, for example, \cite{Schneider:Pn-rk3} for explicit examples).
Thus Theorem \ref{thm:P3} shows that such slope-semistable sheaves 
become unstable with respect to tilt-stability.

In fact, we suggest an even stronger inequality:
\begin{Con} \label{con:strong-BG}
Any tilt-stable object $E \in \BB_{\omega, B}$ with
$\nu_{\omega, B}(E) = 0$ satisfies
\begin{equation} 	\label{eq:strong-BG}
\cht_3(E) \le \frac{\omega^2}{18} \cht_1(E).
\end{equation}
\end{Con}
Just as in the case of the classical Bogomolov-Gieseker inequality for slope-stability,
we have equality when $\omega$ and $B$ are scalar
multiples of the class of an ample line bundle $L$, and $E = L^{\otimes n}$ is a tensor power
of $L$.
We prove this conjecture in the following situations:
\begin{description}
\item[Section \ref{subsec:P3}] For any complex on $\P^3$ when $B = 0$ and
$\omega^3 < 3\sqrt{3}$.
\item[Section \ref{subsec:TorsionSheaves}]
Restrictions of torsion-free sheaves to an ample divisor proportional to $\omega$.
\item[Section \ref{subsec:LineBundles}]
Slope-semistable vector bundles $\FF$ with vanishing discriminant
$\Delta(\FF) = 0$ and $c_1(\FF)$ proportional to $\omega$. (Such sheaves are also stable
with respect to tilt-stability.)
\item[Example \ref{ex:L-otimes-IC}] Sheaves of the form $\OO_X(1) \otimes I_C$ for a curve $C$ on
a hypersurface $X \subset \P^4$, in which case the inequality is related to Castelnuovo's
classical bound for the genus of curves of fixed degree.
\end{description}

If true, the inequality \eqref{eq:strong-BG} would be quite strong.
For example, it would give strong Hodge type inequalities for tilt-stable line bundles when
the N\'eron-Severi group has rank $>1$.
Moreover, in a forthcoming paper \cite{BBMT:Fujita} we show that Conjecture~\ref{con:strong-BG}
implies a Reider-type theorem for threefolds, and a statement towards
Fujita's conjecture on very ampleness of adjoint line bundles (including, for example, Fujita's
conjecture for Calabi-Yau threefolds). The approach is based on ideas and questions in
\cite{AB:Reider}.

As a first step towards proving an inequality for $\ch_3$ in the general case, it seems worthwhile
to generalize the classical Bogomolov-Gieseker for $\ch_0, \ch_1, \ch_2$ from sheaves to complexes.
Indeed, it is an ingredient in any proof of inequalities for $\ch_3$ of slope-stable sheaves; see
\cite{Langer:Survey} for a survey of such inequalities.
Theorem \ref{Thm:Bog} and its Corollaries give various forms of inequalities for tilt-stable
complexes similar to Bogomolov-Gieseker; for example:

\begin{repCor}{Cor:Bog3}
Suppose that $X$ is a smooth projective threefold with N\'eron-Severi group $\NS(X)$ of rank one.
Then any tilt-semistable object $E\in \BB_{\omega, B}$ satisfies
\begin{align*}
\omega (\ch_1(E)^2 -2\rk(E)\ch_2(E)) \ge 0.
\end{align*}
\end{repCor}

Stability conditions at the large-volume limit had been previously constructed in~\cite{large-volume}
and \cite{Toda:limit-stable} as ``polynomial'' or ``limit'' stability condition. As an additional
confirmation that the heart $\AA_{\omega, B}$ seems to give the right construction, we prove:
\begin{Prop}[{Lemma \ref{Lem:L1} and Lemma \ref{Lem:L2}}]
The limit of $\AA_{m\omega, B}$ for $m \to +\infty$ agrees with the heart of polynomial or limit
stability conditions at the large-volume limit.
\end{Prop}
We also prove a compatibility of stability for large $m$ and stability at the limit,
see Proposition \ref{prop:large-volume}.

\subsection{Relation to existing work}\label{subsec:work}
Our construction of $\BB_{\omega, B}$ is directly adopted from Bridgeland's construction
of stability conditions on K3 surfaces in \cite{Bridgeland:K3}. To prove that $\nu_{\omega, B}$
defines a slope function, we use the Bogomolov-Gieseker inequality and the Hodge Index Theorem
just as in the case of general projective surfaces treated by Arcara, Bertram and Lieblich in
\cite{Aaron-Daniele}.
Our notion of tilt-stability on $\BB_{\omega, B}$ is very similar to the notion of
``$\mu_{s + it}$-stability'' by Arcara and Bertram, see \cite{AB:Reider}.

In Section \ref{subsec:P3}, we rely on the construction of ``algebraic stability conditions'' for 
varieties with complete exceptional collections (cf.~\cite{Macri:stability-examples}).
However, even in the case of $\P^3$, our construction includes stability conditions that are not
algebraic.

For large $\omega$, our conjectural stability conditions $(Z_{\omega, B}, \AA_{\omega, B})$ 
should live in a neighborhood of the large volume limit.  Evidently, our approach
is motivated by the string theory construction of $\pi$-stability at the large-volume
limit, see e.g.~ \cite{Douglas:stability, Aspinwall-Douglas:stability,
Aspinwall-Lawrence:DC-zero-brane, Aspinwall:Dbranes-CY}.  In particular, our central charge
\eqref{3fold-charge} is borrowed from the mathematical physics literature \cite[equation
(2.9)]{Aspinwall-Douglas:stability}, with the modification that our formula omits quantum
corrections and a factor of $\sqrt{\td X}$.  This change partly is motivated by the surface case, where
one obtains stability conditions for \emph{every} ample class $\omega$ in this way.

We also refer to \cite{DRY:attractors}
for a conjectural approach to \emph{sufficient} rather than \emph{necessary}
Bogomolov-Gieseker type inequalities on Calabi-Yau threefolds.


\subsection{Acknowledgements} 
We would like to thank Dan Abramovich, Tom Bridgeland, Tommaso de Fernex, Daniel Huybrechts, Mart\'i Lahoz, Adrian Langer, Alexander Polishchuk, Hokuto Uehara for useful discussions and comments.
We are especially grateful to Aaron Bertram and Gueorgui Todorov for many discussions on
Bogomolov-Gieseker type inequalities. We also would like to thank the referee for the detailed
reading of the manuscript.
A.~B.~ is partially supported by NSF grants DMS-0801356/DMS-1001056 and DMS-1101377.  He is also grateful to the
Isaac Newton Institute and its program on ``Moduli Spaces'', during which this paper was finished.
E.~M.~ is partially supported by NSF grant DMS-1001482/DMS-1160466, and by the Hausdorff Center for Mathematics, Bonn, and SFB/TR 45.
Y.~T.~is supported by World Premier International Research Center Initiative (WPI initiative), MEXT, Japan, and Grant-in AId for Scientific Research grant (22684002), partly (S-19104002), from the Ministry of Education, Culture, Sports, Science and Technology, Japan.

\subsection{Notation and Convention}
We work over the complex numbers.
For a set of objects $S$ in a triangulated category $\DD$, we denote by $\langle S \rangle$ the
additive category generated by $S$ via extensions.
If $X$ is smooth and projective variety, we will denote by $\D$ the local dualizing functor
on its derived category $\Db(X)$ given by 
\[
\D(\blank) := (\blank)^\vee[1] := \RlHom(\blank, \OO_X[1]).
\]
Given a coherent sheaf $\FF$, we write $\dim \FF$ for the dimension of its support.
We write $\Coh^{\le d} X =\stv{\FF}{\dim \FF \le d} \subset \Coh X$ for the subcategory of sheaves supported in dimension 
$\le d$, and $\Coh^{\ge d+1} X \subset \Coh X$ for the subcategory of sheaves that have no subsheaf
supported in dimension $\le d$.
Given a bounded t-structure on $\Db(X)$ with heart $\AA$ and an object $E\in\Db(X)$, we write $H^j_{\AA}(E)$, $j\in\mathbb{Z}$, for the cohomology objects with respect to $\AA$.
When $\AA=\Coh(X)$, we simply write $H^j(E)$.

For a complex number $z\in\mathbb{C}$, we denote its real and imaginary part by $\Re z$ and $\Im z$,
respectively. We write $m \gg 0$ to mean ``for \emph{all} sufficiently large $m$''.

We write $\Num(X)$ for the group of cycles $A(X)$ up to numerical equivalence, and $\NS(X) = \NS(X, \Z) = \Num^1(X)$ for the N\'eron-Severi group of divisors up to numerical equivalence.
We also write $\Num_\Q(X)$, $\NS_\Q(X)$, $\Num_\R(X)$, etc.~ for $\Num(X)\otimes\Q$, etc.

We will use the terms ``slope-stability'' and ``$\mu_{\omega, B}$-stability'', as well as
``tilt-stability'' and ``$\nu_{\omega, B}$-stability'' interchangeably when the choice of
$\omega, B$ is clear in context.

\section{Background on stability conditions}\label{sec:Stability}

\subsection{Motivation}
The notion of stability condition on triangulated categories has been introduced by
Bridgeland~\cite{Bridgeland:Stab}, motivated by Douglas's work on
$\Pi$-stability~\cite{Douglas:stability}. We briefly recall the definition:

\begin{Def} A \emph{(full numerical) stability condition} on $\Db(X)$ is a pair $(Z, \AA)$, where 
$\AA \subset \Db(X)$ is the heart of a bounded t-structure and $Z \colon K(\Db(X)) \to \C$ a group
homomorphism, satisfying the following properties:
\begin{enumerate}
\item $Z$ satisfies the positivity property of equation \eqref{eq:Z-positivity}.
\item For the induced notion of stability on $\AA$, every non-zero $E \in \AA$ has a
Harder-Narasimhan filtration in semistable objects in $\AA$.
\item $Z$ factors via the Chern character $\ch \colon K(\Db(X)) \to \Num_\Q(X)$.
\item $(\AA, Z)$ satisfies the ``support property''.
\end{enumerate}
\end{Def}

The support property property will be discussed in the next section.
The set $\Stab(X)$ of stability conditions is a finite-dimensional complex manifold.
In the case where $X$ is a Calabi-Yau threefold, it is expected to contain the stringy
K\"ahler moduli space. More precisely, it should contain an $\Aut \Db(X)$-invariant subspace 
$\NN$ that is isomorphic to the Teichm\"uller space of complex structures on the mirror $\widehat X$; the quotient
$\NN/\Aut \Db(X)$ gets identified with the moduli space of complex structures on $\widehat X$.
Thus the notion of stability conditions on $\Db(X)$ adds a very geometric picture to
Kontsevich's homological mirror symmetry~\cite{KontsICM94}.

The space $\Stab(X)$ has been explicitly studied in several situations.  For instance,
see~\cite{Bridgeland:Stab,Macri:curves,Okada:P1} for $\dim X=1$, and
\cite{Bridgeland:K3,HMS:generic_K3s} for K3 surfaces.  The space $\Stab(X)$ can also
described when $X$ is a local Calabi-Yau variety, e.g.~ the total space of the canonical bundle of a
surface.  For instance, see~\cite{Bridgeland:stab-CY,localP2,Ishii-Ueda-Uehara,Toda:stab-crepant_res,Toda:CY-fibrations}.
The case of non-projective complex tori has been studied in \cite{Sven:generic_tori}.

However, there is no known example of a stability condition on a projective Calabi-Yau threefold,
nor any candidate $(Z, \AA)$ expected to be a stability condition.  One of the issues is that we
have few methods or ideas to construct hearts of a bounded t-structures $\AA \subset \Db(X)$ 
for which the positivity condition \eqref{eq:Z-positivity} could be satisfied.

In principle, one should expect Bridgeland stability conditions for any dimension. Consider the
following central charge, where $X$ is a projective variety of any dimension:

\begin{align} \label{eq:nfold-charge}
Z_{\omega, B}(E) &=-\int_{X}e^{-i\omega}\cht(E) \\\nonumber
&=		\sum_{j \ge 0}\frac{(-1)^{j+1}}{(2j)!}\omega^{2j}\cht_{d-2j}(E)
 + i \left(	\sum_{j\ge 0}\frac{(-1)^j}{(2j+1)!} \omega^{2j+1}\cht_{d-2j-1}(E)\right),
\end{align}

\begin{Con}
There exists a heart $\AA_{\omega, B} \subset \Db(X)$ of a bounded t-structure such that
$(Z_{\omega, B}, \AA_{\omega, B})$ are a Bridgeland stability condition.
\end{Con}

One could apply our method iteratively to construct $\AA_{\omega, B}$ as a $(n-1)$-fold tilt of $\Coh X$.
However, this would involve proving a Bogomolov-Gieseker type inequality at every step.

\subsection{Support property} \label{subsec:support-property}
We require our stability conditions $(Z, \AA)$ to satisfy the following additional technical condition:
\begin{description*}
\item[Support property]
There is a constant $C>0$ such that for any $Z$-semistable object $E\in \AA$, we have
\begin{align*}
\lVert \ch(E) \rVert \le C \lvert Z(E) \rvert,
\end{align*}
where $\lVert \ast \rVert$ is a fixed norm on $\Num_\R(X)$.
\end{description*}
The support property for numerical stability conditions is equivalent (cf.~\cite[Proposition B.4]{localP2})
to the notion of a ``full'' stability condition introduced in \cite{Bridgeland:K3};
in particular:
\begin{Thm}[\cite{Bridgeland:Stab, Kontsevich-Soibelman:stability}] \label{thm:Bridgeland-deform}
There is a natural topology on $\Stab(X)$ such that  the map
\begin{align*}
\Stab(X) \to \Hom(\Num_\Q(X), \C), \qquad (Z, \AA)\mapsto Z, 
\end{align*}
is a local homeomorphism.
\end{Thm}

The support property is also essential to ensure that there is a well-behaved wall-crossing phenomenon 
for stability of objects under deformation of the stability condition:
\begin{Prop} \label{prop:HN-chambers}
Given $E \in \Db(X)$, the set of $(Z, \AA) \in \Stab(X)$ for which $E$ is $(Z, \AA)$-stable is
an open subset of $\Stab(X)$. Further, there exists a chamber decomposition of $\Stab(X)$ by
a locally finite set of walls, such that in the open part of every chamber, the Harder-Narasimhan
filtration of $E$ is constant.
\end{Prop}
This statement is proved by the methods of \cite[Section 9]{Bridgeland:K3}.

\subsection{Tilting}\label{subsec:Tilt}
Our strategy for the construction of $\AA_{\omega, B}$ on threefolds is to take a double tilt 
starting from $\Coh X$ or the category of perverse coherent sheaves on $X$.

\begin{Def}[\cite{Happel-al:tilting}]
Let $\AA$ be the heart of a bounded t-structure on a triangulated category $\DD$.
A pair of subcategories $(\TT, \FF)$ is called a \textit{torsion pair} if the following conditions hold:
\begin{enumerate}
\item \label{tp:Hom-vanishing} For any $T\in \TT$ and $F \in \FF$, we have $\Hom(T, F)=0$.
\item \label{tp:ses}
For any $E\in \AA$, there is an exact sequence $0 \to T \to E \to F \to 0$ in $\AA$, with $T \in \TT$ and $F \in \FF$. 
\end{enumerate}
\end{Def}

Given a torsion pair $(\TT, \FF)$ as above, its \textit{tilt} $\AA^{\dag}$ is the subcategory of $\DD$ defined by
\begin{align*}
\AA^{\dag} =\langle \FF[1], \TT \rangle \subset \DD.
\end{align*}
If $\DD = \Db(\AA)$, then we can identify $\AA^\dag$ with the subcategory of two-term complexes
$E^{-1}\xrightarrow{d} E^0$ with $\ker d \in \FF$ and $\cok d \in \TT$.
The following statements are all well-known:

\begin{Prop} \label{prop:tilt-properties}
\begin{enumerate*}
\item 
The category $\AA^{\dag}$ defined above is the heart of a bounded t-structure on $\DD$.
\item 
Whenever the heart $\AA^{\dag}$ of a bounded t-structure on $\DD$ satisfies
$\AA^{\dag} \subset \langle \AA, \AA[1] \rangle$, then
$\TT = \AA \cap \AA^{\dag}, \FF = \AA \cap \AA^{\dag}[-1]$ is a torsion pair in $\AA$, and
$\AA^{\dag}$ is obtained as the tilt of $\AA$  at $(\TT, \FF)$.
\item
Given two torsion pairs $(\TT_1, \FF_1)$ and $(\TT_2, \FF_2)$ in $\AA$, denote
the corresponding tilts by $\AA_1^{\dag}$ and $\AA_2^{\dag}$, respectively.
If $\TT_2 \subset \TT_1$, then there is a torsion pair
$\TT = \langle \FF_1[1], \TT_2 \rangle, \FF = \FF_2 \cap \TT_1$ in
$\AA_1^\dag$, and
$\AA_2^\dag$ is obtained as the tilt of $\AA_1^\dag$ at this torsion pair.
\end{enumerate*}
\end{Prop}
The first statement is \cite[Proposition~2.1]{Happel-al:tilting}. For the second, see e.g.~ 
\cite[Lemma 1.1.2]{Polishchuk:families-of-t-structures}, and the third statement follows directly
from the second.

\section{First construction}\label{sec:Yukinobu}

Let $X$ be a smooth projective threefold over $\mathbb{C}$.
In this section, we give the first construction of the heart $\AA_{\omega, B}$ of a bounded t-structure on $\Db(X)$ as a double tilt starting from $\Coh X$, 
and state our main conjecture.

\subsection{Tilt of $\Coh X$}\label{susec:FirstTilt}
First we start with the case of arbitrary dimension. 
Let $X$ be a $n$-dimensional 
smooth projective variety
over $\mathbb{C}$, and 
take $B \in \NS_\Q (X)$ and an
ample class $\omega \in \NS_\Q (X)$. 
We use the twisted Chern character $\cht =\ch \cdot e^{-B}$.
Notice that, in particular, we have the following explicit expressions:
\begin{align*}
&\cht_0=\ch_0=\mathrm{rk}\\
&\cht_1=\ch_1-B\ch_0\\
&\cht_2=\ch_2-B\ch_1+\frac{B^2}{2}\ch_0.
\end{align*}
The twisted slope $\mu_{\omega, B}$ on $\Coh X$ is defined as follows.
If $E \in \Coh X$ is a torsion sheaf, we set $\mu_{\omega, B}(E)=+\infty$.
Otherwise we set
\begin{align*}
\mu_{\omega, B}(E)=\frac{\omega^{n-1} \cht_1(E)}{\cht_0(E)}.
\end{align*}
The above slope function satisfies the weak see-saw property, 
i.e., for any exact sequence $0 \to F \to E \to E/F \to 0$
in $\Coh X$
with $F, E/F \neq 0$, one of the following conditions holds, 
\begin{align*}
\mu_{\omega, B}(F) \le \mu_{\omega, B}(E) \le \mu_{\omega, B}(E/F), \\
\mu_{\omega, B}(F) \ge \mu_{\omega, B}(E) \ge \mu_{\omega, B}(E/F). 
\end{align*}
(To prove this, observe that if $\cht_0(F) = 0$, then $F$ is a torsion
sheaf with $\omega^{n-1}\cht_1(F) = \omega^{n-1}(\ch(F)) \ge 0$, and similarly
for $E/F$.)

We define $\mu_{\omega, B}$-stability on $\Coh X$ in the following way: 
$E\in \Coh X$ is $\mu_{\omega, B}$-(semi)stable if, for any $0\neq F \subsetneq E$, we have 
\begin{align}\label{ineq:mu}
\mu_{\omega, B}(F) <(\le) \mu_{\omega, B}(E/F). 
\end{align}
\begin{Rem}
Classically, 
$E \in \Coh X$ is defined to be $\mu_{\omega, B}$-(semi)stable 
if $E$ is torsion free and we have 
the inequality $\mu_{\omega, B}(F)<(\le)\mu_{\omega, B}(E)$
for any subsheaf $0\neq F \subsetneq E$ with $E/F$ torsion free. 
Our definition coincides with the classical definition 
if $E$ has positive rank. An 
inequality similar to (\ref{ineq:mu}) is used 
to define weak stability conditions in~\cite[Section~2]{Toda:PTDT}. 
\end{Rem}
It is well-known that 
the $\mu_{\omega, B}$-stability has the Harder-Narasimhan property, i.e., there is a filtration
\begin{align*}
0=E_0 \subset E_1 \subset \cdots \subset E_N=E,
\end{align*}
such that each $F_i=E_i/E_{i-1}$ is $\mu_{\omega, B}$-semistable with $\mu_{\omega,
B}(F_i)>\mu_{\omega, B}(F_{i+1})$ for all $i$.
We set
\begin{align*}
\mu_{\omega, B; \rm{min}}(E)&=\mu_{\omega, B}(F_N), \\
\mu_{\omega, B; \rm{max}}(E)&=\mu_{\omega, B}(F_1).
\end{align*}
Let $(\TT_{\omega, B}, \FF_{\omega, B})$ be the torsion pair on $\Coh X$ defined by
\begin{align*}
\TT_{\omega, B}&=\left\{ E\in\Coh X\,:\, \mu_{\omega, B; \rm{min}}(E)>0\right\}\\
\FF_{\omega, B}&=\left\{ E\in\Coh X\,:\, \mu_{\omega, B; \rm{max}}(E)\le 0\right\}.
\end{align*}

\begin{Def}
We define the abelian category $\BB_{\omega, B}$ to be the tilt of $\Coh X$ with respect to $(\TT_{\omega, B},
\FF_{\omega, B})$, namely
\begin{equation*}
\BB_{\omega, B}=\langle \FF_{\omega, B}[1], \TT_{\omega, B}\rangle.
\end{equation*}
\end{Def}

Let $Z_{\omega, B}$ be a stability function given by (\ref{eq:nfold-charge}). 
By~\cite{Aaron-Daniele}, we have the following result (the case of K3 surfaces was proved earlier in~\cite{Bridgeland:K3}).

\begin{Prop}[\cite{Aaron-Daniele, Bridgeland:K3}]\label{Rem:BB}
If $\dim X=2$, then $(Z_{\omega, B}, \BB_{\omega, B})$ is a stability condition on $\Db(X)$.
\end{Prop}

The key fact in the proof of the above proposition is the following constraint on numerical classes of
slope-semistable sheaves, known as \textit{Bogomolov-Gieseker inequality}
(see \cite{Reid:Bog, Bogomolov:Ineq, Gieseker:Bog} and \cite[Section 3.4]{HL:Moduli}).
\begin{Thm}[Bogomolov, Gieseker] \label{thm:BG}
Let $X$ be a $n$-dimensional smooth projective variety over $\mathbb{C}$ and let $\omega$ be an ample divisor on $X$.
For any torsion free $\mu_{\omega, B}$-semistable sheaf $E$, we have the following inequality:
\begin{align*}
\omega^{n-2}\bigl(\cht_1(E)^2 -2\cht_0(E) \cht_2(E)\bigr) \ge 0. 
\end{align*}
\end{Thm}

\subsection{Tilt of $\BB_{\omega, B}$}\label{susec:SecondTilt}
From now on, we focus on the case $\dim X=3$; as stated in equation \eqref{3fold-charge},
the central charge is then given by
\begin{align*}
Z_{\omega, B}(E)=\left(-\cht_3(E)+\frac{\omega^2}{2}\cht_1(E) \right)
+i\left(\omega \cht_2(E)-\frac{\omega^3}{6}\cht_0(E)  \right).
\end{align*} 
The abelian category $\BB_{\omega, B}$ satisfies the following property:

\begin{Lem} \label{lem:def-nu}
For any non-zero object $E\in \BB_{\omega, B}$, one of the following conditions holds:
\begin{enumerate}
\item $\omega^2 \cht_1(E)  >0$.
\item $\omega^2 \cht_1(E) =0$ and $\Im Z_{\omega, B}(E)>0$.
\item $\omega^2 \cht_1(E)= \Im Z_{\omega, B}(E)=0$ and $-\Re Z_{\omega, B}(E)>0$.
\end{enumerate}
\end{Lem}

\begin{Prf}
By the construction of $\BB_{\omega, B}$, we have $\omega^2 \cht_1(E) \ge 0$.
Suppose that $\omega^2 \cht_1(E)=0$.
Then $H^0(E) \in \Coh^{\le 1} X$ and $H^{-1}(E)$ is
 $\mu_{\omega, B}$-semistable torsion free sheaf with $\mu_{\omega, B}(E)=0$.
By the Hodge Index Theorem and the Bogomolov-Gieseker inequality, we have
\begin{align*}
0\ge
 \omega \cht_1(H^{-1}(E))^2 \ge 2 \omega \cht_0(H^{-1}(E)) \cht_2(H^{-1}(E)),
\end{align*}
which implies $\omega \cht_2(H^{-1}(E)) \le 0$.
Since $\cht_0(E) \le 0$ and $\omega \cht_2(H^0(E)) \ge 0$, we obtain the inequality $\Im Z_{\omega, B}(E) \ge 0$.
Finally suppose that $\omega^2 \cht_1(E)=\Im Z_{\omega, B}(E)=0$.
Then the above argument shows that $H^{-1}(E) = 0$ and $E = H^0(E)$ has zero-dimensional support; hence the inequality $-\Re Z_{\omega, B}(E)>0$ holds.
\end{Prf}

\begin{Rem}
The above lemma implies that the vector $(\omega^2 \cht_1, \Im Z_{\omega, B}, -\Re Z_{\omega, B})$ 
for objects of $\BB_{\omega, B}$ behaves like the vector $(\ch_0, \ch_1, \ch_2)$ for coherent sheaves on a surface.
The subcategory of $E\in \BB_{\omega, B}$ satisfying $\omega^2 \cht_1(E)=0$ is an analogue of the
subcategory of torsion sheaves; we can also describe it as the extension-closure
\begin{align}\label{026}
\langle \Coh^{\le 1} X, F[1]\, \text{for all $\mu_{\omega, B}$-stable $F$ with } \mu_{\omega, B}(F)=0 \rangle.
\end{align}
In case $B = 0$, the above category contains the heart of the category of D0-D2-D6 bound states constructed in~\cite{Toda:PTDT}.
\end{Rem}

We define a slope $\nu_{\omega, B}$ on $\BB_{\omega, B}$ as follows.
If $E\in \BB_{\omega, B}$ satisfies $\omega^2 \cht_1(E)=0$, we set $\nu_{\omega, B}(E)=+\infty$.
Otherwise we set
\begin{align} \label{eq:nu-def}
\nu_{\omega, B}(E) := \frac{\Im Z_{\omega, B}(E)}{\omega^2 \cht_1(E)} 
= \frac{\omega \cht_2(E) - \frac 16 \omega^3 \cht_0(E)}{\omega^2 \cht_1(E)}.
\end{align}
By Lemma~\ref{lem:def-nu}, 
the slope $\nu_{\omega, B}$ also  
satisfies the weak see-saw property. 
Therefore an analogue of slope stability on $\BB_{\omega, B}$
is defined in the following way:
\begin{Def}
An object $E\in \BB_{\omega, B}$ is $\nu_{\omega, B}$-(semi)stable if, for any non-zero proper
subobject $F\subset E$ in $\BB_{\omega, B}$, we have the inequality
\begin{align*}
\nu_{\omega, B}(F)<(\le) \mu_{\omega, B}(E/F). 
\end{align*}
\end{Def}

Similarly to $\mu_{\omega, B}$-stability, we have the following result.

\begin{Lem}
The Harder-Narasimhan property holds with respect to $\nu_{\omega, B}$-stability, i.e., for any
$E\in \BB_{\omega, B}$, there is a filtration in $\BB_{\omega, B}$
\begin{align}\label{HN:v}
0=E_0 \subset E_1 \subset \cdots \subset E_N=E,
\end{align}
such that each $F_i=E_i/E_{i-1}$ is $\nu_{\omega, B}$-semistable with $\nu_{\omega,
B}(F_i)>\nu_{\omega, B}(F_{i+1})$ for all $i$.
\end{Lem}

\begin{Prf}
First we note that $\BB_{\omega, B}$ is a noetherian abelian category.
This is essentially proved in~\cite{Bridgeland:K3} when $X$ is a K3 surface, and almost the same proof works in the general case.
Indeed, we only need to modify the argument 
of~\cite[Prop. 7.1]{Bridgeland:K3} in the following way. 
In the notation of~\cite[Prop. 7.1]{Bridgeland:K3}, 
the sheaves $H^0(L_i)$ turned out to be the finite length 
sheaves in the K3 surface case. In our 3-fold situation, 
the sheaves $H^0(L_i)$ are at most one dimensional, 
so may not be
of 
finite length. 
However, since the codimensions 
of the supports of $H^0(L_i)$ are at least two, 
we obtain a chain
\begin{align*}
H^{-1}(E_1) \subset H^{-1}(E_2) \subset \cdots \subset Q^{\ast \ast},
\end{align*}
in the notation of~\cite[Prop. 7.1]{Bridgeland:K3}.
Instead of bounding the length of $H^0(L_i)$, we 
can terminate the above chain as $\Coh(X)$ is noetherian.
This proves that $\BB_{\omega, B}$ is noetherian. 

Since $B$ and $\omega$ are rational, we can then apply the same arguments as in~\cite[Prop.\ B.2]{localP2} to show the Harder-Narasimhan property. 
\end{Prf}

For an object $E\in \BB_{\omega, B}$ 
with Harder-Narasimhan filtration (\ref{HN:v}) we set
\begin{align*}
\nu_{\omega, B; \rm{min}}(E)&=\nu_{\omega, B}(F_N), \\
\nu_{\omega, B; \rm{max}}(E)&=\nu_{\omega, B}(F_1),
\end{align*}
and the torsion pair $(\TT'_{\omega, B}, \FF'_{\omega, B})$ on $\BB_{\omega, B}$ is defined by
\begin{align*}
\TT_{\omega, B}'&=\left\{ E\in\BB_{\omega, B} \,:\, \nu_{\omega, B; \rm{min}}(E)>0\right\}\\
\FF_{\omega, B}'&=\left\{ E\in\BB_{\omega, B} \,:\, \nu_{\omega, B; \rm{max}}(E)\le 0\right\}.
\end{align*}

\begin{Def}
We define the abelian category $\AA_{\omega, B}$ to be the tilt of $\BB_{\omega, B}$
with respect to $(\TT_{\omega, B}', \FF_{\omega, B}')$, namely
\begin{align*}
\AA_{\omega, B}=\langle \FF_{\omega, B}'[1], \TT_{\omega, B}' \rangle.
\end{align*}
\end{Def}

By the construction of $\AA_{\omega, B}$, it is obvious that $\Im Z_{\omega, B}(E) \ge 0$, for all $E\in \AA_{\omega, B}$.
We propose the following conjecture. 

\begin{Con}\label{Con:stability}
The pair $(Z_{\omega, B}, \AA_{\omega, B})$ is a stability condition on $\Db(X)$.
\end{Con}

The above conjecture in particular implies that, for any $\nu_{\omega, B}$-semistable object $E\in
\BB_{\omega, B}$ with $\nu_{\omega, B}(E)=0$, we have $\Re Z_{\omega, B}(E)>0$.
More precisely, Conjecture~\ref{Con:stability} immediately implies the following conjecture.

\begin{Con}\label{Con:stability2}
For any $\nu_{\omega, B}$-semistable object $E\in \BB_{\omega, B}$ satisfying
\begin{align*}
\frac{\omega^3}{6}\cht_0(E) =\omega \cht_2(E),
\end{align*}
we have
\begin{align*}
\cht_3(E)<\frac{\omega^2}{2}\cht_1(E).
\end{align*}
\end{Con}

In Section \ref{sec:Noether} we will show that the two conjectures are equivalent,
by showing that $\AA_{\omega, B}$ is Noetherian.

\subsection{Support property for tilt-stability}
To show that there is a well-behaved notion of wall-crossing for tilt-stability of objects
$E \in \Db(X)$, we need some form of boundedness of potentially destabilizing subobjects. This
boundedness follows from a form of the ``support property'' discussed in Section
\ref{subsec:support-property}. To set this up, define a central charge
$\overline{Z}_{\omega, B} \colon K(\Db(X)) \to \C$
corresponding to the slope function $\nu_{\omega, B}$:
\[
\overline{Z}_{\omega, B}(E) = \frac 12 \omega^2 \cht_1(E) + i \Im Z_{\omega, B}(E).
\]

\begin{Rem} \label{rem:Zbar}
By \eqref{026}, we have
$\overline{Z}_{\omega, B}(E) \in \stv{r e^{i\pi\phi}}{r > 0, -\frac 12 < \phi \le \frac 12}$
for every non-zero $E \in \BB_{\omega, B}$ except if $E$ is a zero-dimensional torsion sheaf.
As, for such objects, the slope induced by $\overline{Z}_{\omega,B}$ agrees with $\nu_{\omega, B}$, an object
$k(x)\neq E \in \BB_{\omega, B}$ is tilt-stable if and only if $\Hom(k(x), E) = 0$ and there are no 
destabilizing subobjects with respect to $\overline{Z}_{\omega, B}$.
\end{Rem}

\begin{Lem}
Fix a norm $\lVert \cdot \rVert$ on $\Num_\Q(X)$. 
There exists a constant $C > 0$ such that, for every tilt-stable object $E \in \Db(X), 
E \neq k(x)[n]$, we have
\begin{equation}\label{eq:SupportProperty}
\lVert \ch(E) \rVert \le C \lvert \overline{Z}_{\omega,B}(E) \rvert
\end{equation}
\end{Lem}
\begin{Prf}
We give a sketch of the argument;
the complete proof is in \cite[Sections 3.6 \& 3.7]{Toda:ExtremalContractions}.
Using the same methods as in the proof of the support property for surface, given in \cite[Section
4]{localP2}, we will show \eqref{eq:SupportProperty} only for the semi-norm $\lVert \ch
\rVert':=\lVert (\ch_0,\omega^2\ch_1,\omega\ch_2,\ch_3) \rVert$, which will 
be enough for all applications (in particular in section \ref{sec:BGwithoutch3}). The full 
statement can then be deduced from Theorem \ref{Thm:Bog}. 

For any torsion-free slope-stable sheaf $\FF$, define $x_{\omega, B}(\FF), y_{\omega, B}(\FF)\in\R$ by
\[
x_{\omega, B} + iy_{\omega, B} = \frac{\overline{Z}_{\omega, B}(\FF)}{\rk \FF}.
\]
Using the classical Bogomolov-Gieseker inequality and the Hodge inequality as in the proof of
Lemma \ref{lem:def-nu}, one shows that
\begin{align}\label{eq:Columbus16}
y_{\omega, B} = -\frac {\omega^3}6 + \frac{\omega \cht_2(\FF)}{\rk(\FF)}
\le -\frac {\omega^3}6 + \frac{\omega \cht_1(\FF)^2}{2 \rk(\FF)^2}
\le -\frac {\omega^3}6 + \frac{2x^2}{\omega^3} =: f_\omega(x)
\end{align}

We define a function $S_{\min}$ of $\omega,B$ by
\[
S_{\min}(\omega, B) = \inf \stv{\abs{x + i f_\omega(x)}}{x \in \R}
\]
The function $S_{\min}$ is continuous, and by \eqref{eq:Columbus16} it satisfies
\[
0 < S_{\min}(\omega,B) \leq \inf \left\{ \left|\frac{\overline{Z}_{\omega,B}(\FF)}{r(\FF)} -i t\right|\,:\,
			t \in \R_{\ge 0},\, \FF\text{ torsion-free slope-stable sheaf}\right\}.
\]
Using exactly the same arguments as in \cite[Lemma 4.5]{localP2} one then deduces
the claim for objects where either $H^0(E)$ or $H^{-1}(E)$ have positive rank.
Using the openness of the ample cone, the claim also follows for torsion sheaves.
\end{Prf}

We could also formulate tilt-stability completely in the formalism of weak stability conditions
introduced in \cite{Toda:PTDT}. Then the support property would be satisfied for every stable
object, including $k(x)$. Since we are not interested in deforming the slope of skyscraper sheaves
of points, the above Lemma is sufficient for our purposes:

\begin{Cor} \label{cor:tilt-wall-crossing}
Denote by $U \subset \NS_\R(X) \times \NS_\R(X)$ the set of pairs $(\omega, B)$ where $\omega$ is ample.
The notion of tilt-stability can be extended to all pairs $(\omega, B) \in U$.
For every object, the set of $(\omega, B) \in U$ for which $E$ is tilt-stable is an open subset 
of $U$. Further, there is a chamber decomposition of $U$, given by a locally finite set of walls, such that
the Harder-Narasimhan filtration of $E$ is constant on every chamber.
\end{Cor}
\begin{Prf}
The first claim follows from Bridgeland's deformation result recalled in 
Theorem \ref{thm:Bridgeland-deform}. As in Proposition \ref{prop:HN-chambers}, it follows that there
is a chamber decomposition for stability with respect to $\overline{Z}_{\omega, B}$. Combined
with Remark \ref{rem:Zbar}, this implies the claim.
\end{Prf}

\section{Second construction}\label{sec:polynomial}

The second construction of the heart $\AA_{\omega, B}$ starts from perverse
coherent sheaves rather than sheaves, and uses polynomial stability conditions rather than
slope-stability. We will compare the two constructions in Section \ref{sec:Comparison}.

\subsection{Polynomial stability conditions}\label{subsec:prelim}
The notion of polynomial stability condition has been introduced in \cite{large-volume}. We refer to \emph{loc. cit.} for all basic definitions.
We will repeatedly construct polynomial stability conditions by using the following proposition/definition - which is stated slightly differently in \cite{large-volume}, but the proof is the same.

\begin{PropDef}\label{def:polynomial}
Let $\DD$ be a triangulated category.
Giving a polynomial stability condition on $\DD$ is equivalent to giving a heart of
a bounded t-structure $\AA \subset \DD$, and a central charge
$Z \colon K(\DD) \to \C[m]$ such that
\begin{enumerate}
\item \label{cond-stabfunction}
For every $0 \neq E \in \AA$, and for some fixed $a \in \R$, the leading coefficient
of $Z(E)$ is contained in the semi-closed half plane
\[
\R_{>0} \cdot e^{i \pi (a, a+1]}.
\]
\item \label{enum:HN-exists}
Harder-Narasimhan filtrations exist for the stability condition on $\AA$ induced by $Z$.
\end{enumerate}
\end{PropDef}

We say that $Z$ is a \emph{stability function} with respect to the interval $(a, a+1]$
if it satisfies condition \eqref{cond-stabfunction}. In this case, we can define a ``polynomial
phase function'' for every $E \in \AA$: it is a continuous function germ
$\phi(E) \colon (\R \cup {+\infty}, +\infty) \to \R$ defined by 
\[
\phi(E)(m) = \frac 1{\pi} \arg Z(E)(m)
\]
for sufficiently large $m$, such that
\begin{align*}
\lim_{m\to +\infty} \phi(E)(m) \in (a, a+1].
\end{align*}
 For two such functions $\phi, \phi'$ we say
$\phi \succ \phi'$ if $\phi(m) > \phi'(m)$ for $m \gg 0$. This defines a notion of stability for
objects in $\AA$, by comparing its polynomial phase function with that of its subobjects,
and condition \eqref{enum:HN-exists} of Definition \ref{def:polynomial} refers to HN-filtrations
with respect to this notion of stability.

For a polynomial stability condition $(Z, \PP)$, the \emph{slicing} $\PP$ gives the set of semistable objects
for every polynomial phase functions $\phi$.
We let $\widehat{\PP}$ be the induced $\R$-valued slicing given by
\[
\widehat{\PP}(\phi) = \langle \stv{\PP(\phi(m))}{\phi(+\infty) = \phi} \rangle.
\]
More concretely, in the setting of Proposition \ref{def:polynomial}, and for $\phi \in (a, a+1]$,
the subcategory $\widehat{\PP}(\phi) \subset \AA$ is extension-closure generated by $Z$-semistable
objects $E$ with 
\begin{align*}
\lim_{m \to +\infty} \frac 1{\pi} \arg Z(E)(m) = \phi.
\end{align*}
The key input of polynomial stability conditions is that, having constructed a polynomial $(Z, \PP)$
from a heart $\AA$ as above, we get new t-structures by setting $\AA' := \OPP((b, b+1])$ for
any $b \in \R$. The category $\AA'$ could also be described as (the shift of) a tilt of
$\AA$.

We will repeatedly use the following lemma, which is established in the proof of
\cite[Theorem 2.29]{Toda:limit-stable}. We refer to \cite[Section 4]{Bridgeland:Stab} for the notion
of quasi-abelian categories.
\begin{Lem}\label{lem:torsionpair-HN}
Let $\TT, \FF$ be a torsion pair in $\AA$, and $Z$ a polynomial stability function for $\AA$.
Write $\phi(E)$ for the polynomial phase functions induced by $Z$ on $\AA$.
Assume that
\begin{enumerate}
\item \label{enum:compareslopes}
If $T \in \TT, F \in \FF$, then $\phi(T) \succ \phi(F)$.
\item HN-filtrations exist for $Z$ on the quasi-abelian categories $\TT, \FF$.
\end{enumerate}
Then HN-filtrations exist for $Z$ on $\AA$.
\end{Lem}
In this situation, an object $E \in \AA$ is $Z$-stable if and only if
$E \in \TT$ or $E \in \FF$, and it is $Z$-stable in the respective quasi-abelian category with
respect to strict inclusions.

Finally, we recall the notion of ``dual stability condition'': We say that the polynomial stability
conditions $(Z', \PP')$ and $(Z, \PP)$ are dual to each other if $\PP'(\phi) = \D(\PP(-\phi))$, and
if $Z'(\D(E))$ is the complex conjugate of $Z(E)$. Recall that we included a shift by one in
the definition $\D(E) = E^\vee[1]$ of our dualizing functor.
Hence if skyscrapher sheaves
$k(x)$ are stable with respect to $\PP$ of phase 1, then the corresponding statement holds for
$\PP'$: 
\[ \D(k(x)) = k(x)[-2] \in \PP(-1) = \PP'(1). \]

\subsection{Perverse stability}\label{sec:perverse}
The starting point is a polynomial stability condition on the category of perverse coherent sheaves: 
\begin{Def}\label{def:perverse}
We define the category of perverse coherent sheaves $\Coh^p$ to be  
\begin{align*}
\Coh^p =\langle \Coh^{\ge 2} X [1], \Coh^{\le 1} X \rangle,
\end{align*}
and the central charge
 $Z_p^{\omega, B} \colon K(X) \to \mathbb{C}[m]$ to be  
\begin{equation} \label{eq:def-Zp}
Z_p^{\omega, B}(E) := 
 - \cht_3(E) + m i \omega \cht_2(E) + 
	m^2 \left(\frac{\omega^2}{2} \cht_1(E) - i \frac{\omega^3}{6} \cht_0(E) \right).
\end{equation}
\end{Def}

\begin{figure}
\begin{center}
        \subfloat[$(Z_p, \Coh^p)$]{ \includegraphics{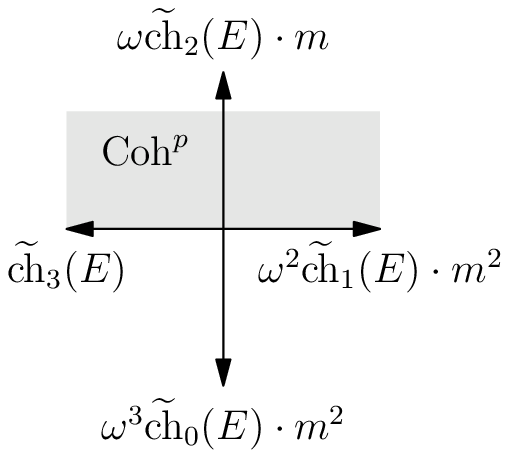} }
        \subfloat[$(Z_\BB, \BB)$]{ \includegraphics{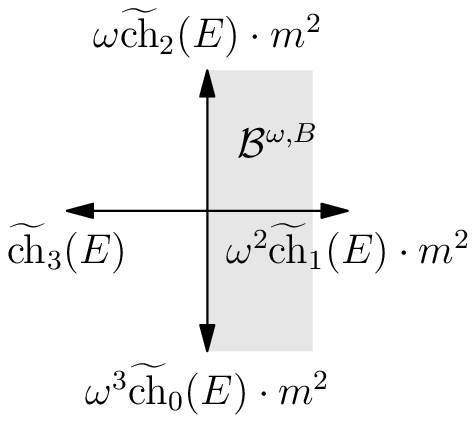} }
\caption{The auxiliary polynomial stability conditions}
\label{fig:Zp-ZB}
\end{center}
\end{figure}

Our strategy, also indicated in Figure~\ref{fig:Zp-ZB}, is as follows: first we show that
$(Z_p^{\omega, B}, \Coh^p)$ gives a polynomial stability condition whose heart corresponds to the
upper half plane. The central charge is dominated by the $\cht_0$ and $\cht_1$-terms; in other
words, this stability condition is a refinement of slope-stability.
In the next step, we rescale the contribution of $\cht_2$ to have the same 
weight of $m^2$; as this only changes the imaginary part of the central charge, this is done after
switching to the tilt $\BB^{\omega, B}$ of $\Coh^p$, which is the heart corresponding to the right
half-plane. The resulting stability is closely related to tilt-stability.

In the final (and conjectural step), we rescale the contribution of $\cht_3$;
since that only changes the real part of the central charge, we first switch to the tilt
$\AA^{\omega, B}$ that corresponds to the upper half plane.

To show $(Z_p^{\omega, B}, \Coh^p)$ is a stability condition, we follow the method of \cite[Theorem
2.29]{Toda:limit-stable}; more precisely, we will use the following results about a torsion
pair in $\Coh^p$:
\begin{Lem}[{\cite[Lemma~2.16, Lemma~2.17 and Lemma~2.19]{Toda:limit-stable}}] \label{lem:Yukinobu-torsionpair}
There exists a torsion pair 
$(\AA_{1}^p, \AA_{1/2}^p)$ in $\Coh^p$ defined as follows:
\begin{align*}
\AA^p_{1} &= \langle F[1], k(x) \,:\, F \mbox{ is pure two-dimensional, } x\in X \rangle, \\
\AA^p_{1/2} &= \left\{E \in \Coh^p  \,:\, \Hom(\AA_{1}^p, E)=0  \right\}.
\end{align*}
Each of the quasi-abelian categories $\AA^p_1$ and $\AA^p_{1/2}$ is of finite length with respect to
strict inclusions and strict epimorphisms.

Additionally, they satisfy $\D(\AA^p_1) = \AA^p_1[-2]$ and $\D(\AA^p_{1/2}) = \AA^p_{1/2}[-1]$.
\end{Lem}

\begin{Prop}\label{prop:HNZp}
$Z_p^{\omega, B}$ is a stability function for $\Coh^p$, and HN-filtrations exist.
\end{Prop}
\begin{Prf}
If $H^{-1}(E)$ for an object $E \in \AA^p_{1/2}$ does not vanish,
then $H^{-1}(E)$ is purely three-dimensional; hence $E$ has negative rank, and the leading term of
$Z_p^{\omega, B}(E)$ has positive imaginary part. Similarly, if $H^{-1}(E) = 0$, then
$H^0(E)$ is purely one-dimensional, and the same conclusion holds.
On the other hand, for an object $E \in \AA^p_{1}$ the leading term evidently has a negative real
coefficient.

This shows that $Z_p^{\omega, B}$ is a stability function on $\Coh^p$
with respect to the interval $(0, 1]$, and at the same time that $(\AA^p_1, \AA^p_{1/2})$ satisfies the condition
\eqref{enum:compareslopes} in Lemma \ref{lem:torsionpair-HN}. Since $\AA^p_1$ and $\AA^p_{1/2}$ are
of finite length, the existence of HN-filtrations for each of them is also satisfied, and 
the conclusion follows from Lemma \ref{lem:torsionpair-HN}.
\end{Prf}

Let $\PP_p^{\omega, B}$ be the induced slicing with values in ``polynomial phase functions'', and
$\widehat{\PP}_p^{\omega, B}$ the induced $\R$-valued slicing as defined above.

\begin{Prop} \label{prop:Zp-selfdual}
The stability condition $(Z_p^{\omega, B}, \PP_p^{\omega, B})$ is dual to the stability condition
$(Z_p^{\omega, -B}, \PP_p^{\omega, -B})$.
\end{Prop}
\begin{Prf}
The Chern characters of $E$ and $\D(E)$ differ by a sign in $\ch_0$ and $\ch_2$, and agree for
$\ch_1$ and $\ch_3$; the same holds for 
$\ch(E) e^B$ and $\ch(\D(E)) e^{-B}$. Thus 
\begin{equation} \label{eq:Zpduality}
 \overline{Z_p^{\omega, -B}(\D(E))} = Z_p^{-\omega, -B}(\D(E)) = 
Z_p^{\omega, B}(E).
\end{equation}

Furthermore, in the proof of Proposition \ref{prop:HNZp}, we identified $\AA^p_1$ with
$\widehat{\PP}_p^{\omega, B}(1)$ and $\AA^p_{1/2}$ with $\widehat{\PP}_p^{\omega, B}((0, 1))$. Combined with the last
statement of Lemma \ref{lem:Yukinobu-torsionpair}, this implies that
\begin{align*}
\D(\widehat{\PP}_p^{\omega, B}(1)) = \widehat\PP_p^{\omega, -B}(-1), \
\D(\widehat{\PP}_p^{\omega, B}((0, 1)) = 
\widehat\PP_p^{\omega, -B}((-1, 0)).
\end{align*}
As the $\D$ turns
strict inclusions in these quasi-abelian categories into strict epimorphisms, and vice versa,
equation \eqref{eq:Zpduality} also implies that $E$ is $Z_p^{\omega, B}$-stable if and only if $\D(E)$ is
$Z_p^{\omega, -B}$-stable.
\end{Prf}

\subsection{Surface-like stability}\label{sec:surface}
We use polynomial stability condition $(Z_p^{\omega, B}, \PP_p^{\omega, B})$ to define the heart
$\BB^{\omega, B}$.

\begin{Def}
Let $\BB^{\omega, B} := \widehat{\PP}_p^{\omega, B}((-\frac 12, \frac 12])$ and
\[ 
Z_\BB^{\omega, B} := 
 - \cht_3(E) + m^2 \left(\frac 12 \omega^2 \cht_1(E) + i \left(\omega \cht_2(E) -  \frac 16 \omega^3 \cht_0(E)\right) \right).
\]
\end{Def}
In other words, $\BB^{\omega, B}[1]$ is the tilt of $\Coh^p$ at the torsion pair
$\widehat{\PP}_p^{\omega, B}((\frac 12, 1]), \widehat{\PP}_p^{\omega, B}((0, \frac 12])$.

\begin{Prop} 			\label{prop:ZB-stability}
$Z_\BB^{\omega, B}$ is a stability function for $\BB^{\omega, B}$ with respect to the interval $(-\frac 12, \frac 12]$ and HN-filtrations exist.
\end{Prop}

To prove Proposition \ref{prop:ZB-stability} we need a more detailed understanding of the cohomology
sheaves for objects in $\BB^{\omega, B}$ (a more precise result will be Lemma \ref{lem:Comparison1}):

\begin{Lem}\label{lem:Cohomologies1}
The cohomology sheaves of any $Z_p^{\omega, B}$-stable object $E \in \BB^{\omega, B}$ either vanish, or satisfy
the following conditions:
\begin{enumerate}
\item\label{cohom:1} $H^1(E)$ is a zero-dimensional torsion sheaf.
\item\label{cohom:2} $H^{-1}(E)$ is a slope-stable torsion-free sheaf of slope $\mu_{\omega, B} \leq 0$.
\item\label{cohom:3} $H^0(E)$ is either a slope-semistable torsion-free sheaf of slope $\mu_{\omega, B} > 0$, or a torsion sheaf. Moreover, if $H^0(E)$ has a zero-dimensional subsheaf, then $H^{-1}(E)$ is also non-zero.
\end{enumerate}
\end{Lem}

\begin{Prf}
This follows from the following statements about $Z_p^{\omega, B}$-stable objects
$E \in \Coh^p$:
\begin{itemize}
\item If $H^0(E)$ has one-dimensional support, then $E \in 
\OPP_p^{\omega, B}((0, \frac 12])$. 
\item If $H^{-1}(E) \neq 0$, then either $E \in \OPP_p^{\omega, B}(1)$ and $H^{-1}(E)$ is a purely two-dimensional
sheaf, or $H^{-1}(E)$ is slope-semistable. Its slope $\mu_{\omega, B}$ satisfies $\mu_{\omega, B} \le 0$
if and only if $E \in \OPP_p^{\omega, B}((0, \frac 12])$. 
\end{itemize}
Indeed, there is a surjection $E \onto H^0(E)$ in $\Coh^p$, which destabilizes $E$ unless the first claim
holds. To show the second claim, first assume that $H^{-1}(E) \in \Coh^{\ge 2}$ is purely
three-dimensional but not slope-semistable, and
let $A \subset H^{-1}(E)$ be a destabilizing subsheaf. Then the composition $A[1] \into
H^{-1}(E)[1] \into E$ is an inclusion in $\Coh^p$ that destabilizes $E$.
The same argument deals with the case where $H^{-1}(E)$ is not purely three-dimensional.
This shows directly \eqref{cohom:1} and \eqref{cohom:2}.
To prove \eqref{cohom:3}, we only need to observe that, if $H^0(E)$ has a torsion subsheaf of dimension zero, then this destabilizes $E$ unless $H^{-1}(E)$ is also non-zero.
\end{Prf}

\begin{Prf} (Proposition \ref{prop:ZB-stability})
To prove that $Z_\BB^{\omega, B}$ is a stability function first note that
$\Re Z_\BB^{\omega, B} = \Re Z_p^{\omega, B}$, and that if the leading coefficient of $Z_p^{\omega, B}(E)$ has positive real part, then
the same holds for $Z_\BB^{\omega, B}(E)$. In particular, if $E \in \OPP_p^{\omega, B}((-\frac 12, \frac 12))$, then 
$Z_p^{\omega, B}(E)$ has leading coefficient with positive real part, and so $Z_\BB^{\omega, B}(E)$ satisfies
the required property.

In the remaining case we have $E \in \OPP_p^{\omega, B}(\frac 12)$. If $H^{-1}(E) \neq 0$, then it is a
slope-semistable sheaf of slope $\mu_{\omega, B} = 0$. From the Bogomolov-Gieseker inequality it follows that
$\omega \cdot \cht_2(H^{-1}(E)) \le 0$. Additionally, for any $E \in \Coh^p$ with $\cht_0(E)=\cht_1(E)=0$, we have
that $\omega \cdot \cht_2(E) \ge 0$. It follows that
the leading coefficient of $Z_\BB^{\omega, B}(E)$ is a positive imaginary number.

To prove the existence of HN-filtrations, first note that the torsion pair given by $\TT = \OPP_p^{\omega, B}(\frac 12), \FF = \OPP_p^{\omega, B}((-\frac 12, \frac 12))$ satisfies
condition \eqref{enum:compareslopes} of Lemma \ref{lem:torsionpair-HN}. Due to the rationality of
$\omega$ and of $t$, the imaginary part of $Z_\BB^{\omega, B}$ is discrete, and thus
$\OPP_p^{\omega, B}(\frac 12)$ has finite length.
By the following Lemma, 
the quasi-abelian category $\FF = \OPP_p^{\omega, B}((-\frac 12, \frac 12))$ is also of finite
length, and thus our claim follows from Lemma \ref{lem:torsionpair-HN}.
\end{Prf}

\begin{Lem} \label{lem:FFfinite-length}
The quasi-abelian category $\OPP_p^{\omega, B}((-\frac 12, \frac 12))$ has finite length.
\end{Lem}
\begin{Prf}
As above, we denote this category by $\FF$.
By Proposition \ref{prop:Zp-selfdual}, the dual $\D(\FF)$ is of the same form as $\FF$ itself
(with $B$ replaced by $-B$); thus, it is enough to check that there are no infinite chains
$\dots \into E_3 \into E_2 \into E_1 $ of strict subobjects in $\FF$. By the rationality of
$\omega$ and $s$, we may assume that the real part of the $m^2$-coefficient of
$Z_\BB^{\omega, B}(E_j)$ is constant. But
then the imaginary part of the $m^2$-coefficient must also be constant, as the quotient $Q_j$ of
$E_{j+1} \into E_j$ could otherwise not lie in $\OPP_p^{\omega, B}((-\frac 12, \frac 12))$. In particular,
$Z_\BB^{\omega, B}(Q_j)$ is a constant polynomial.

From the proof above of the fact that $Z_\BB^{\omega, B}$ is a stability function it follows that this is only
possible if $Z_p^{\omega, B}(Q_j)$ already was a constant polynomial, which means that $Q_j$ is the shift
$T[-1]$ of a zero-dimensional skyscraper sheaf. Hence the long exact cohomology sequence 
induces a sequence of inclusions $H^1(E_{j+1}) \into H^1(E_j)$ of zero-dimensional torsion sheaves,
which must terminate.
\end{Prf}

For later use, we also show a partial converse to Lemma \ref{lem:Cohomologies1}:

\begin{Lem}\label{lem:skyscr}
If $T$ is a torsion sheaf of dimension zero, then $T[-1]\in\BB^{\omega,B}$.
Moreover, if $E\in\BB^{\omega,B}$, then the exact triangle
\begin{equation}\label{eq:Bonn07022011}
Q\to E\to H^1(E)[-1]
\end{equation}
gives an exact sequence in $\BB^{\omega,B}$, where $Q$ is the extension
\[
H^{-1}(E)[1]\to Q \to H^0(E).
\]
\end{Lem}
\begin{Prf}
Let $T$ be a torsion sheaf of dimension zero.
Then $T[-1]$ belongs to $\BB^{\omega,B}$ since it is stable of phase $0$ with respect to $Z_p^{\omega, B}$.

By Lemma \ref{lem:Cohomologies1}, by looking at the long exact sequence for the cohomology sheaves, an object $M\in\BB^{\omega,B}$ is a subobject of $T[-1]$ if and only if $M[1]$ is a torsion sheaf of dimension zero and $M[1]\into T$ in $\Coh(X)$: for an exact sequence in $\BB^{\omega,B}$
\[
0\to M\to T[-1]\to N\to0,
\]
we have
\[
0\to H^0(N)\to H^1(M)\to T\to H^1(N)\to 0,
\]
and $H^{-1}(N)\cong H^0(M)$, which is impossible unless
\[
H^{-1}(N)=H^0(N)=H^{-1}(M)=H^0(M)=0.
\]
To show that \eqref{eq:Bonn07022011} gives an exact sequence in $\BB^{\omega,B}$ it is enough to observe that, if the non-zero map $E\to H^1(E)[-1]$ is not surjective in $\BB^{\omega,B}$, then, by what we just proved, it must factorize through a torsion subsheaf of $H^1(E)$, which is clearly a contradiction.
\end{Prf}

We write $(Z_\BB^{\omega, B}, \PP_\BB^{\omega, B})$ for the induced polynomial stability condition,
and $\widehat\PP_\BB^{\omega, B}$ for the corresponding $\R$-valued slicing.

\begin{Prop}\label{prop:ZB-selfdual}
The stability condition $(Z_\BB^{\omega, B}, \PP_\BB^{\omega, B})$ is dual to the stability
condition $(Z_\BB^{\omega, -B}, \PP_\BB^{\omega, -B})$.
\end{Prop}
\begin{Prf}
Observe that, by the construction of $\BB^{\omega, B}$, we have 
\begin{align*}
\OPP_\BB^{\omega, B} \left(\pm \frac{1}{2} \right)
&=\OPP_{p}^{\omega, B} \left(\pm \frac{1}{2} \right), \\ 
\OPP_\BB^{\omega, B} \left(\left(-\frac{1}{2}, \frac{1}{2} \right)\right)
&=\OPP_{p}^{\omega, B} \left(\left(-\frac{1}{2}, \frac{1}{2} \right) \right).\\
\intertext{Also by Lemma~\ref{prop:Zp-selfdual}, we have}
\D\left(\widehat\PP_p^{\omega, B}\left(\frac 12\right)\right) &= 
\widehat\PP_p^{\omega, B}\left(-\frac 12 \right), \\
\D\left(\widehat\PP_p^{\omega, B}\left(\left(-\frac 12, \frac 12\right)\right)\right) &= 
\widehat\PP_p^{\omega, B}\left(\left(-\frac 12, \frac 12\right)\right).
\end{align*}
Then the claim follows with the same arguments as Proposition \ref{prop:Zp-selfdual}.
\end{Prf}

\subsection{The threefold heart}
\begin{Def} We define
$\AA^{\omega, B}$ to be the heart 
\[ \AA^{\omega, B}:=\widehat{\PP}_\BB^{\omega, B}((0,1]) \]
of the slicing $\widehat{\PP}_\BB^{\omega, B}$. 
\end{Def}
Note that since $\BB^{\omega, B} = \widehat{\PP}_\BB^{\omega, B}((-\frac 12, \frac 12])$, one could also define
$\AA^{\omega, B}$ as the tilt of $\BB^{\omega, B}$ at the torsion pair 
$\TT = \widehat{\PP}_\BB^{\omega, B}((0, \frac 12])$, 
$\FF = \widehat{\PP}_\BB^{\omega, B}((-\frac 12, 0])$.

Let $Z_{\omega, B}$ be the central charge
defined by (\ref{3fold-charge}). 
By construction, the imaginary part of $Z_{\omega, B}$ and $\frac 1{m^2} Z_\BB^{\omega, B}$ agree;
thus we automatically have 
\begin{align*}
\Im Z_{\omega, B}(E) \ge 0,
\ E \in \AA^{\omega, B}.
\end{align*}
 To show that $Z_{\omega, B}$ is a
stability function on $\AA^{\omega, B}$, we would have to show that objects $E \in \AA^{\omega, B}$ with 
$\Im Z_{\omega, B}(E) = 0$ satisfy $\Re Z_{\omega, B}(E) < 0$; equivalently, if $E \in \BB^{\omega, B}$ is 
$Z_\BB^{\omega, B}$-stable of phase 0, then $\Re Z_{\omega, B}(E) > 0$. 
In the next section, we will prove  that $\AA^{\omega, B}$ equals $\AA_{\omega, B}$, and
so this claim is equivalent to Conjecture \ref{Con:stability2}.

\begin{Rem}\label{rmk:cohomologies}
Let $E \in \BB^{\omega, B}$ be a $Z_\BB^{\omega, B}$-stable object of phase 0, which is not isomorphic to $k(x)[-1]$. Then $E$ is quasi-isomorphic to a two-term complex $E^{-1} \to E^0$ of vector bundles. Indeed, as $k(x)[-1]$ is
stable of the same phase, we have $\Hom(E, k(x)[n]) = 0$ for $n \le -1$ and (using Serre duality)
for $n \ge 2$; then the claim follows by \cite[Corollary 5.6]{Bridgeland-Maciocia:K3fibrations}.
\end{Rem}

\begin{Rem}
If the stability condition $(\AA^{\omega, B}, Z^{\omega, B})$ exists, then it is
dual to $(\AA^{\omega, -B}, Z^{\omega, -B})$; the proof is the same as for
Proposition \ref{prop:ZB-selfdual}. 
\end{Rem}

\section{Comparison and Noetherian property}

\subsection{Comparing the two constructions}\label{sec:Comparison}
The goal of this section is to prove the following result:

\begin{Prop}\label{prop:Comparison}
We have $\AA_{\omega, B}=\AA^{\omega, B}$.
\end{Prop}

The proof consists of a detailed analysis of the various steps in the two constructions.

\medskip

\noindent
{\bf Step 1.} ($\Coh^p$ versus $\Coh X$)
By definition, $\Coh^p$ is the tilt of $\Coh X$ with respect to the torsion pair
\begin{align*}
\TT_0&=\Coh^{\leq 1} X\\
\FF_0&=\Coh^{\geq 2} X.
\end{align*}

\medskip

\noindent
{\bf Step 2.} ($\BB^{\omega, B}$ versus $\BB_{\omega, B}$)
By Lemma \ref{lem:Cohomologies1} and Lemma \ref{lem:skyscr}, the skyscraper sheaves $k(x)[-1]\in\BB^{\omega, B}$ and, for all $E\in\BB^{\omega, B}$, $H^1(E)$ is torsion of dimension $0$.
Define a torsion pair in $\BB^{\omega, B}$ by
\begin{align*}
\TT_1&=\left\{ E\in\BB^{\omega, B}\,:\,H^1(E)=0\right\}\\
\FF_1&=\left\{ E\in\BB^{\omega, B}\,:\,E\cong H^1(E)[-1]\right\}.
\end{align*}
Notice that, by its own definition, $\FF_1[1]=\Coh^{\leq 0}X$ consists of zero-dimensional sheaves.
The fact that this is a torsion pair follows immediately from Lemma \ref{lem:skyscr}.
Let $\BB_1$ be the tilt with respect to this torsion pair, i.e.,
\[
\BB_1=\langle \FF_1[1], \TT_1 \rangle.
\]

\begin{Lem}\label{lem:Comparison1}
We have $\BB_1=\BB_{\omega, B}$.
\end{Lem}
\begin{Prf}
We only need to show that $\BB_{\omega, B}\subseteq\BB_1$.
Let $M\in\BB_{\omega, B}$.
By construction of $\BB_{\omega, B}$, both of its cohomology sheaves $H^0(M)$ and $H^{-1}(M)[1]$ belong
to $\BB_{\omega, B}$. By further using their Harder-Narasimhan filtrations with respect to
$\mu_{\omega, B}$-stability, it is sufficient to consider the following cases:

\begin{enumerate}
\item \label{enum:torsion}
$M = \Gamma$ is a torsion sheaf.
\item \label{enum:bigslope}
$M = \Gamma$ is a torsion-free slope-stable sheaf with slope $\mu_{\omega, B}(\Gamma)>0$,
\item \label{enum:smallslope}
$M = \Gamma[1]$ is the shift of a torsion-free slope-stable sheaf $\Gamma$ with slope 
$\mu_{\omega, B}(\Gamma) \le 0$,
\end{enumerate}

For case (\ref{enum:torsion}), we can assume $\Gamma$ is pure and so we can distinguish three sub-cases,
according to the dimension of the support of $\Gamma$:
\begin{description*}
\item[$\dim(\Gamma)=0$] In this case, $\Gamma[-1] \in \FF \subset \BB^{\omega, B}$ and $\Gamma\in\BB_1$ by
construction of $\BB_1$.
\item[$\dim(\Gamma)=1$] The limit phase with respect to $Z_p^{\omega, B}$ is $\frac 12$.
Assume, for a contradiction, that $\Gamma\notin
\OPP_p^{\omega, B}((0,\frac 12])$.
Then there exists an exact sequence in $\Coh^p$
\[
A\to \Gamma\to B
\]
where $A$ is $Z_p^{\omega, B}$-semistable with limit phase $\overline{\phi}_p(A)>\frac 12$.
Passing to cohomology, we have
\[
0\to H^{-1}(B)\to H^0(A)\to \Gamma \to H^0(B)\to 0.
\]
Since $\dim(H^{-1}(B))\geq2$ and $\dim(H^0(A)) \le 1$, we have $H^{-1}(B)=0$.
Hence $A\cong H^0(A)$ must be pure of dimension $1$ and its limit phase $\overline{\phi}_p(A)$ is precisely $\frac 12$, a contradiction.
\item[$\dim(\Gamma)=2$] The limit phase of $\Gamma[1] \in \Coh^p$ with respect to $Z_p^{\omega, B}$ is $1$.
If $A \into \Gamma[1] \onto B$ is an exact sequence in $\Coh^p$, then the long exact cohomology
sequence shows
that $B[-1]$ is a sheaf, and in fact that it is a sheaf with 2-dimensional support. Thus any
quotient of $\Gamma[1]$ also has limit phase $1$, and so $\Gamma[1]\in\OPP^{\omega, B}_p(1)$.
But then $\Gamma\in\BB^{\omega, B}\cap\Coh X\subset\BB_1$.
\end{description*}

In case (\ref{enum:bigslope}) we have $\Gamma[1]\in\Coh^p = \OPP_p^{\omega, B}((0, 1])$.
Assume that $\Gamma[1]\notin\OPP_p^{\omega, B}((\frac 12,1])$.
Then there exists an exact sequence in $\Coh^p$
\[
A \to \Gamma[1] \to B
\]
where $B$ is $Z_p^{\omega, B}$-semistable and has limit phase $\overline{\phi}_p(B)\leq\frac 12$.
Passing to cohomology, we have
\[
0\to H^{-1}(A)\to \Gamma\xrightarrow{\phi} H^{-1}(B)\to H^0(A)\to 0.
\]
Hence $B\cong H^{-1}(B)[1]$. By Lemma \ref{lem:Cohomologies1},
$H^{-1}(B)$ is torsion-free and $\mu_{\omega, B}$-semistable with slope $\mu_{\omega, B}(H^{-1}(B))\leq0$.
But then $\phi=0$ and $H^{-1}(B)\cong H^0(A)=0$.

This contradiction proves that $\Gamma[1]\in\OPP_p^{\omega, B}((\frac 12,1])$.
Hence, $\Gamma \in \OPP_p^{\omega, B}((-\frac 12, 0])\subset \BB^{\omega, B}$, and clearly
it is contained in $\TT_1 \subset \BB_1$.

Finally we treat case (\ref{enum:smallslope}). Consider the exact sequence in $\Coh^p$
\begin{equation}\label{eqn:Tucson1}
A\to\Gamma[1]\to B
\end{equation}
with $A \in \OPP_p^{\omega, B}
((\frac 12, 1])$ and $B \in \OPP_p^{\omega, B}((0, \frac 12])$.
Passing to cohomology, we have
\[
0\to H^{-1}(A)\to \Gamma\xrightarrow{\phi} H^{-1}(B)\to H^0(A)\to 0.
\]
If $H^{-1}(A)$ is non-zero, then its slope satisfies $\mu_{\omega, B}(H^{-1}(A)) \leq
\mu_{\omega, B}(\Gamma)\leq0$. 
Since $H^0(A)$ is a torsion sheaf of dimension $\leq1$, we have
$\overline{\phi}_p(A) = \overline{\phi}_p(H^{-1}(A))\le \frac 12$, which is a contradiction; hence
$A = H^0(A)$.
If $H^0(A)$ has dimension $1$, its limit phase is $\overline{\phi}_p(H^0(A))=\frac 12$, which is
again a contradiction.
Hence $T_0:=H^0(A)$ is a 0-dimensional torsion sheaf.
Thus the exact sequence \eqref{eqn:Tucson1} becomes
\[
T_0\to\Gamma[1]\to\Lambda[1],
\]
where $\Lambda$ is a sheaf; in particular, $\Lambda[1] \in \OPP_p^{\omega, B}((0, \frac 12]) \subset
\BB^{\omega, B}$ 
is contained in the torsion-part of $\BB^{\omega, B}$ and thus $\Lambda[1] \in \BB_1$. 
We already proved that $T_0\in\BB_1$ in part (\ref{enum:torsion}), and thus
we also have $\Gamma[1] \in \BB_1$.
\end{Prf}

\medskip

\noindent
{\bf Step 3.} ($\AA^{\omega, B}$ versus $\AA_{\omega, B}$)
It will be enough to show that $\AA_{\omega, B}\subseteq\AA^{\omega, B}$. Of course, the key point
will be that the slope of the $m^2$-coefficient of $Z_\BB^{\omega, B}$ is, up to normalization,
given by $\nu_{\omega, B}$.

By Step 2 and Lemma \ref{lem:skyscr}, there is a torsion pair $(\TT_2, \FF_2)$ in $\BB_{\omega, B}$
where $\TT_2 = \FF_1[1] = \Coh^{\le 0}X$ consists of zero-dimensional skyscraper sheaves; then
$\FF_2 = \BB_{\omega, B} \cap \BB^{\omega, B}$ is given as the right-orthogonal
\[
\FF_2 = \stv{E \in \BB_{\omega, B}}{\Hom(k(x), E) = 0,\, \text{for all $x \in X$}}
\]
Evidently $\TT_2 \subset \AA^{\omega, B}$. Recall that $\nu_{\omega, B}(k(x)) = +\infty$.
Hence if $E \subset \BB_{\omega, B}$ is tilt-semistable,
then either $\nu_{\omega, B}(E) = +\infty$ and the short exact sequence
$T \into E \onto E'$ with $T \in \TT_2$ and $E' \in \FF_2$ has $E'$ also tilt-semistable
with $\nu_{\omega, B}(E') = +\infty$, or $E$ itself is already in $\FF_2$.
It follows that it is sufficent to show
for every tilt-stable object $E \in \FF_2$:
\begin{enumerate}
\item \label{enum:nu<=0}
If $\nu_{\omega, B}(E) \le 0$, then $E \in \OPP_\BB^{\omega, B}((-\frac 12, 0])$.
\item \label{enum:nu>0}
If $\nu_{\omega, B}(E) > 0$, then $E \in \OPP_\BB^{\omega, B}((0, \frac 12])$.
\end{enumerate}
Indeed, by definition, $\AA^{\omega, B}=\OPP_{\BB}^{\omega, B}((0,1])$, and thus the claim
implies $\TT'_{\omega, B} \subset \AA^{\omega, B}$ and $\FF'_{\omega, B}[1] \subset \AA^{\omega, B}$.

Consider such an $E$. By Step 2, we have $E \in \TT_1 \subset \BB^{\omega, B}$. Consider
any short exact sequence 
\[
A \into E \onto B
\]
in $\BB^{\omega, B}$ that destabilizes $E$ with respect to
$Z_{\BB}^{\omega, B}$. Consider the short exact sequence $T \into A \onto A/T$ given by the 
torsion pair $(\TT_1, \FF_1)$ in $\BB^{\omega, B}$. If $T$ is non-trivial, then the limiting phases
of $Z_\BB^{\omega, B}(T)$ and $Z_\BB^{\omega, B}(A)$ agree (as $Z_\BB^{\omega, B}(A/T)$ is a constant
polynomial).
On the other hand, consider the induced short exact sequence
\[ T \into E \onto E/T. \]
As $T, E \in \TT_1$, and as $\TT_1 \subset \BB^{\omega, B}$ is closed under quotients,
this is also a short exact sequence
in $\BB_{\omega, B}$. By the tilt-stability of $E$ we have $\nu_{\omega, B}(T) < \nu_{\omega,
B}(E)$. This is a contradiction unless $T = 0$.

Hence either $E$ is stable with respect to $Z_\BB^{\omega, B}$, or its Harder-Narasimhan filtration
has just two steps $0 \into E_1 \into E_2$ with
$E_1 \in \FF_1 = \Coh^{\le 0}[-1]$ being the shift of a zero-dimensional skyscraper sheaf.

In case \eqref{enum:nu<=0}, the limiting phase of $Z_\BB(E)(m)$ satisfies
$\phi(E)(+\infty) = \phi(E_2/E_1)(+\infty) \in (-\frac 12, 0]$; together with $\phi(k(x)[-1]) = 0$ this 
shows our claim. In the other case \eqref{enum:nu>0}, we have $\phi(E)(+\infty) \in (0, \frac 12]$.
In particular, $\phi(E) \succ \phi(E_1) = \phi(k(x)[-1]) = 0$, a contradiction unless $E_1 = 0$.
Thus $E$ is $Z_\BB^{\omega, B}$-stable with $E \in \OPP_\BB^{\omega, B}((0, \frac 12])$.

This finishes the proof of Proposition \ref{prop:Comparison}.
We also have the following more precise result, which will be used in \cite{BBMT:Fujita}:

\begin{Prop}
\begin{enumerate}
\item \label{enum:precise-comp}
Assume $E\in\BB_{\omega, B}\cap\BB^{\omega, B}$ satisfies
either $\nu_{\omega, B}(E) \ge 0$ or
$\Hom(k(x)[-1], E) = 0$. Then $E$ is $\nu_{\omega, B}$-semistable
if and only if it is $Z_\BB^{\omega, B}$-semistable.

\item Assume $E \in \BB_{\omega, B}$ satisfies $\nu_{\omega, B; \max}(E) < +\infty$. 
Then $\D(E)$ fits into an exact triangle 
$ \widetilde E \to \D(E) \to T_0[-1]$
for an object $\widetilde E\in \BB_{\omega, -B}$ and a
zero-dimensional torsion sheaf $T_0$.  Further, $E$ is $\nu_{\omega, B}$-semistable if and only if
$\widetilde E$ is $\nu_{\omega, -B}$-semistable.
\end{enumerate}
\end{Prop}
\begin{Prf}
The first claim follows from the proof of Step 3. For the second claim, first note that
$\nu_{\omega, B; \max}(E) < +\infty$ implies
$\Hom(k(x), E) = 0$ and thus $E \in \BB_{\omega, B} \cap \BB^{\omega, B}$. Combined with the first
part, this gives the even stronger statement
$E \in \OPP^{\omega, B}_\BB((-\frac 12, \frac 12))$. Proposition \ref{prop:ZB-selfdual} then implies
$\D(E) \in \OPP^{\omega, -B}_\BB((-\frac 12, \frac 12)) \subset \BB^{\omega, -B}$, and Lemma
\ref{lem:Comparison1} implies the existence of an exact triangle as stated above.

If $E$ is not $\nu_{\omega, B}$-stable, then it has a destabilizing quotient
$E \onto B$ with $\nu_{\omega, B; \max}(B) < +\infty$. Applying the same construction to $B$ produces
an injection $\widetilde B \into \widetilde E$ in $\BB_{\omega, -B}$. As
$\nu_{\omega, B}(E) = - \nu_{\omega, -B}(\widetilde E)$ etc., this will destabilize $\widetilde E$
with respect to $\nu_{\omega, -B}$. So $E$ unstable implies $\widetilde E$ unstable, and the
converse follows similarly.
\end{Prf}

\subsection{Noetherian property} \label{sec:Noether}
The goal of this section is to show that the heart of our t-structure is Noetherian.

\begin{Lem}\label{lem:Access}
There is a torsion pair $(\TT^1, \FF^{(0,1))}$ in $\AA_{\omega, B}$ whose torsion part is given by
\[
E \in \TT^1 \quad \Longleftrightarrow \quad \Im Z_{\omega, B}(E) = 0.
\]
The category $\TT^1$ is an abelian category of finite length which is closed in $\AA_{\omega,B}$ under subobjects and quotients.
\end{Lem}
\begin{Prf}
Using the second construction $\AA^{\omega, B} = \OPP_\BB^{\omega, B}((0, 1])$, we can define a torsion pair
$\TT^1 := \OPP_\BB^{\omega, B}(1)$ and $\FF^{(0,1)} = \OPP_\BB^{\omega, B}((0, 1))$.
Evidently $E \in \TT^1$ if and only if
the leading coefficient of $Z_\BB^{\omega, B}(E)$ is real, which happens if and only if
$Z_{\omega, B}(E)$ is real.

To prove the second assertion, notice that $\TT^1$ is by definition a subcategory of the quasi-abelian category $\OPP_\BB^{\omega, B}((\frac 12, \frac 32))$.
By construction of $(Z_\BB^{\omega, B}, \PP_\BB^{\omega, B})$, we have 
$\OPP_\BB^{\omega, B}((\frac 12, \frac 32)) = \OPP_p^{\omega, B}((\frac 12, \frac 32))$, and the latter is of finite length by Lemma \ref{lem:FFfinite-length}.
\end{Prf}

\begin{Prop} \label{prop:Noether}
The abelian category $\AA_{\omega, B}$ is Noetherian.
\end{Prop}
\begin{Prf}
Suppose that there is an infinite sequence of surjections in $\AA_{\omega, B}$
\begin{align}\label{infi}
E_1 \twoheadrightarrow E_2 \twoheadrightarrow \cdots.
\end{align}
We are going to show that the above sequence terminates.
Since $\omega$ and $B$ are rational and
$\Im Z_{\omega, B}(E) \ge 0$ for any $E\in \AA_{\omega, B}$, we may assume that
\begin{align}\label{Zo}
\Im Z_{\omega, B}(E_1)=\Im Z_{\omega, B}(E_i),
\end{align}
for all $i$. Consider the  exact sequence
\begin{align*}
0 \to L_i \to E_1 \to E_i \to 0
\end{align*}
in $\AA_{\omega, B}$.
By equation \eqref{Zo} we have $\Im Z_{\omega, B}(L_i) = 0$ and so $L_i \in \TT^1$.
Thus every $L_i$ is a subobject
of the torsion part $T$ of $E_1$. Replacing $E_1$ by $T$ and $E_i$ by the quotient
$T/L_i$, we get an infinite sequence \eqref{infi} with $E_i \in \TT^1$.
This sequence terminates by the second part of Lemma \ref{lem:Access}.
\end{Prf}

\begin{Rem}
If $(Z_{\omega, B}, \AA_{\omega, B})$ is a stability condition, then $\AA_{\omega, B}$ must be Noetherian as shown in~\cite[Proposition~10.1]{Abramovich-Polishchuk:t-structures}.
Hence Proposition~\ref{prop:Noether} gives some evidence for Conjecture~\ref{Con:stability}.
\end{Rem}

\begin{Cor}\label{cor:HNfiltrations}
Conjecture~\ref{Con:stability} and Conjecture~\ref{Con:stability2} are equivalent.
\end{Cor}

\begin{Prf}
It is obvious that Conjecture~\ref{Con:stability} implies Conjecture~\ref{Con:stability2}.
Suppose that Conjecture~\ref{Con:stability2} is true.
Then $Z_{\omega, B}$ is a stability function on $\AA_{\omega, B}$.
Since $\AA_{\omega, B}$ is Noetherian by Proposition~\ref{prop:Noether}, we can
apply~\cite[Prop.~B.2]{localP2} to conclude that $(Z_{\omega, B}, \AA_{\omega, B})$ has the Harder-Narasimhan property.
\end{Prf}

\section{Large volume limit}\label{sec:LVL}

In this section we show that if the stability conditions $(\AA_{\omega, B}, Z_{\omega, B})$ exist, then their
limit as $\omega$ goes to infinity is exactly given by the notion of ``polynomial stability
condition at the large-volume limit'' of \cite[Section 4]{large-volume} or ``limit stability'' of
\cite{Toda:limit-stable}. The precise statement is given in Proposition \ref{prop:large-volume}.

\subsection{Stability condition at the large volume limit}\label{subsec:StabilityLVL}
Let $Z_{\infty \omega, B} \colon K(X) \to \mathbb{C}[m]$
be a polynomial valued central charge, given by 
\begin{align*}
Z_{\infty \omega, B}(E)(m)=
Z_{m\omega, B}(E). 
\end{align*}
(Note that the only difference to $Z_p^{\omega, B}$ of equation \eqref{eq:def-Zp}
is given by the $\cht_0(E)$-term, which has weight $m^3$ rather than $m^2$.)
Let $\Coh^p$ be the category of perverse 
coherent sheaves, given in Definition~\ref{def:perverse}. 
Recall that, as in Section \ref{sec:perverse}, by \cite{large-volume} the pair $(Z_{\infty\omega, B}, \Coh^p)$
 defines a polynomial stability condition on $\Db(X)$.
Let $\QQ^{p}$ 
be the associated slicing depending on polynomial phase functions. 
The central charge $Z_{\infty\omega, B}$ is a stability 
function on $\Coh^p$ with respect to the interval
$(1/4, 5/4)$ (see \cite[Lemma~2.20]{Toda:limit-stable}), 
hence we have 
$\Coh^p=\QQ^p ((1/4, 5/4))$.
\begin{Def}
We define $\CC^p:=\QQ^p((0, 1])$.
\end{Def}

We give a precise description of the abelian category $\CC^p$. 
Note that there is an analogue of slope stability on $\Coh^{\le 2} X$.
Namely, for an object $E\in \Coh^{\le 2} X$, we set $\widehat{\mu}_{\omega, B}(E)=+\infty$ if $E\in \Coh^{\le 1} X$, and otherwise we set
\begin{align} \label{eq:def-muhat}
\widehat{\mu}_{\omega, B}(E)=\frac{\omega \cht_2(E)}{\omega^2 \cht_1(E)}.
\end{align}
The $\widehat{\mu}_{\omega, B}$-stability on $\Coh^{\le 2} X$ is defined in a similar way to $\mu_{\omega, B}$-stability. 
We define the torsion pair $(\TT^p, \FF^p)$ on $\Coh(X)$ to be
\begin{align*}
\TT^p&=\left\{ E\in\Coh^{\le 2}
 X\,:\, \widehat{\mu}_{\omega, B; \rm{min}}(E)>0\right\}\\
\FF^p&=\left\{ E\in\Coh X\,:\, \Hom(\TT^p, E)=0\right\}.
\end{align*}

\begin{Lem}
The abelian category $\CC^p$ is the tilt of $\Coh X$ with respect to $(\TT^p, \FF^p)$,
\begin{align*}
\CC^p=\langle \FF^p[1], \TT^p \rangle.
\end{align*}
\end{Lem}

\begin{Prf}
It is enough to show that the RHS is contained is the LHS.
To see this, it is enough to check that
\begin{enumerate}
\item Any $\widehat{\mu}_{\omega, B}$-stable sheaf $E\in \Coh^{\le 2} X$ with
$\widehat{\mu}_{\omega, B}(E)>0$ (resp.~$\widehat{\mu}_{\omega, B}(E) \le 0$) satisfies $E\in \CC^p$ (resp.~$E[1] \in \CC^p$).
\item Any torsion free sheaf $E\in \Coh X$ satisfies $E[1] \in \CC^p$.
\end{enumerate}
Let $(\AA_{1}^p, \AA_{1/2}^p)$ 
be the torsion pair on $\Coh^p$, defined in 
Lemma~\ref{lem:Yukinobu-torsionpair}. 
It is shown in~\cite[Lemma~2.27]{Toda:limit-stable} that an object $E\in \Coh^p$ is $Z_{\infty\omega,
B}$-semistable if and only if $E\in \AA_{i}^p$ for $i=1$ or $1/2$ and it is $Z_{\infty\omega, B}$-semistable in the quasi-abelian category $\AA_i^p$.

Suppose that $E\in \Coh^{\le 2} X$ is $\widehat{\mu}_{\omega, B}$-stable with
$\widehat{\mu}_{\omega, B}(E)>0$. 
If $E$ has two-dimensional support, 
then $E[1] \in \AA^{p}_{1}$ and it is $Z_{\infty\omega, B}$-stable in $\AA^{p}_1$. Since $\Im
Z_{m\omega, B}(E[1]) <0$ for $m\gg 0$, we have $E[1] \in \QQ^p(>1)$, hence $E\in \CC^p$. If $E$ is
pure of dimension one, then $E \in \AA^{p}_{1/2} \subset \CC^p$; and if $E \cong k(x)$
is the skyscraper sheaf of a point $x \in X$, then it is $Z_{\infty\omega, B}$-stable of phase 1.
A similar argument shows that if $\widehat{\mu}_{\omega, B}(E) \le 0$, then $E[1] \in \CC^p$.

Next, take a torsion free sheaf $E\in \Coh X$.
Then $E[1] \in \AA_{1/2}^{p}$, and its Harder-Narasimhan factors with respect to the polynomial
stability function $Z_{\infty\omega, B}$ are contained in $\AA_{1/2}^{p}$.
Since any object in $\AA_{1/2}^{p}$ has limit phase $1/2$, we have $E[1] \in \CC^p$.
\end{Prf}

\begin{figure}
\begin{centering}
\begin{tikzpicture} 
  [node distance=.0cm,
  start chain=going right,]
     \node[punktchain] (A) {$\Coh^{=2}_{\widehat{\mu}_{>0}}[1]$};
     \node[punktchain] (B)      {$\Coh^{=2}_{\widehat{\mu} {\le0}}[1]$};
     \node[punktchain] (C)      {$\Coh^{=3}_{\mu_{> 0}}[1]$};
     \node[punktchain] (D) {$\Coh^{=3}_{\mu_{\le 0}}[1]$};
     \node[punktchain] (E) {$\Coh^{\le 1}$};
     \node[punktchain] (F) {$\Coh^{=2}_{\widehat{\mu}_{>0}}$};
     \node[punktchain] (G)      {$\Coh^{=2}_{\widehat{\mu}_{\le0}}$};
     \node[punktchain] (H)      {$\Coh^{=3}_{\mu_{> 0}}$};
     \node[punktchain] (J) {$\Coh^{=3}_{\mu_{\le 0}}$};
\draw[tuborg] let
    \p1=(D.west), \p2=(H.east) in
    ($(\x1,\y1+1.7em)$) -- ($(\x2,\y2+1.7em)$) node[above, midway]  {$\BB_{\omega, B}$};
\draw[tuborg, decoration={brace}] let
    \p1=(E.west), \p2=(J.east) in
    ($(\x2,\y2-1.7em)$) -- ($(\x1,\y1-1.7em)$) node[below, midway]  {$\Coh X$};
\draw[tuborg, decoration={brace}] let
    \p1=(B.west), \p2=(F.east) in
    ($(\x2,\y2-1.8em)$) -- ($(\x1,\y1-1.8em)$) node[below, midway]  {$\CC^p$};
\draw[tuborg, decoration={brace}] let
    \p1=(A.west), \p2=(E.east) in
    ($(\x1,\y1+1.8em)$) -- ($(\x2,\y2+1.8em)$) node[above, midway]  {$\Coh^p$};
\end{tikzpicture}
\caption{Relating the various t-structures}
\label{fig:tstructures}
\end{centering}
\end{figure}
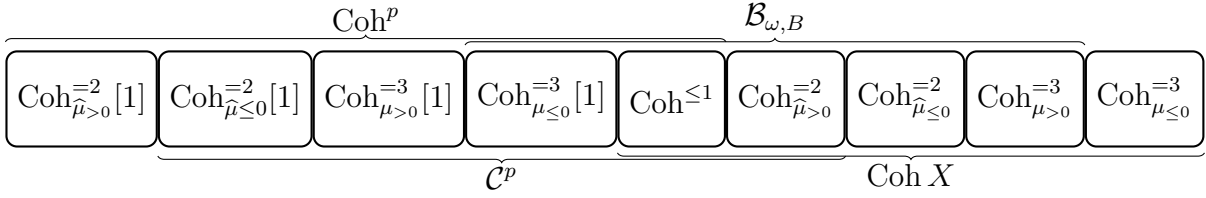

The diagram in Figure~\ref{fig:tstructures} schematically shows the relations between the different 
t-structures. Each heart in the figure is the extension-closure of the corresponding blocks.

\subsection{Comparison of $\AA_{\omega, B}$ and $\CC^p$}\label{subsec:ComparisonLVL}
\begin{Lem}\label{Lem:L1}
For an object $E\in \Db(X)$, suppose that $E\in \AA_{m\omega, B}$ for $m\gg 0$.
Then $E\in \CC^p$.
\end{Lem}

\begin{Prf}
Note that we have $\BB_{\omega, B}=\BB_{m\omega, B}$ for $m \in \mathbb{R}_{>0}$.
We denote by $H_{\BB}^{i}(\blank)$ the $i$-th cohomology functor with respect to the t-structure
$\BB_{\omega, B}$.
Also, for simplicity, we write $\TT_{m\omega, B}'$, $\FF_{m\omega, B}'$ and $\AA_{m\omega, B}$ as
$\TT_{m}', \FF_m'$ and $\AA_{m}$ respectively.
Suppose that $E\in \Db(X)$ satisfies $E\in \AA_m$ for $m \gg 0$.
This implies that
\begin{align}\label{eqn:L1}
H_{\BB}^{-1}(E) \in \FF_m'\qquad \text{ and } \qquad H_{\BB}^{0}(E) \in \TT_m',
\end{align}
for $m \gg 0$.
We have the following exact sequences in $\BB_{\omega, B}$:
\begin{align*}
&0 \to E_1[1] \to H_{\BB}^{-1}(E) \to E_2 \to 0, \\
& 0 \to E_3[1] \to H_{\BB}^{0}(E) \to E_4 \to 0,
\end{align*}
for $E_i \in \Coh X$.
By the construction of $\BB_{\omega, B}$, we have $E_1, E_3 \in \FF_{\omega, B}$ and $E_2, E_4 \in
\TT_{\omega, B}$.
Since $(\TT_{m}', \FF_{m}')$ is a torsion pair on $\BB_{\omega, B}$, (\ref{eqn:L1}) implies that $E_1[1] \in \FF_m'$ and $E_4 \in \TT_m'$ for $m \gg 0$.
In particular we have $\nu_{m\omega, B}(E_1[1])\le 0$ and $\nu_{m\omega, B}(E_4)>0$ for $m\gg 0$, which imply that
\begin{align*}
-\frac{1}{6}m^3 \omega^3 \cht_{0}(E_i)+ m\omega \cht_{2}(E_i) \ge 0,
\end{align*}
for $i=1, 4$ and $m \gg 0$. 
This implies that $\cht_{0}(E_1)=\cht_{0}(E_4)=0$, hence $E_1=0$ and $E_4$ is a torsion sheaf.

From what we have proved above, the object $E$ is concentrated on $[-1, 0]$.
Since $E_3 \in \FF^p$, it is enough to check that
\begin{align*}
E_2 \in \FF^p\qquad \text{ and }\qquad E_4 \in \TT^p,
\end{align*}
to conclude $E\in \CC^p$. 
Let $E_{2, \rm{tor}} \subset E_2$ be the torsion part of $E_2$, and $F\subset E_{2, \rm{tor}}$ be
the $\widehat{\mu}_{\omega, B}$-semistable factor of $E_{2, \rm{tor}}$ with $\widehat{\mu}_{\omega,
B}$ maximum.
Then $F \in \FF_{m}'$ for $m \gg 0$, as $F$ is a subobject of $H_{\BB}^{-1}(E)$ in $\BB_{\omega,
B}$, therefore $\nu_{m\omega, B}(F) \le 0$ for $m \gg 0$.
This implies that $F$ is a pure two-dimensional sheaf with $\widehat{\mu}_{\omega, B}(F) \le 0$, hence $E_2 \in \FF^p$ follows.
Similarly for a $\widehat{\mu}_{\omega, B}$-semistable factor $E_4 \twoheadrightarrow F'$ such that
$\widehat{\mu}_{\omega, B}$ is minimum, we have $F' \in \TT_m'$ for $m \gg 0$, hence $\mu_{\omega,
B}(F')>0$ and $E_4 \in \TT^p$ follows.
\end{Prf}

\begin{Lem}\label{Lem:L2}
For an object $E\in \CC^p$, we have $E\in \AA_{m\omega, B}$ for $m \gg 0$.
\end{Lem}

\begin{Prf}
Let us take an object $E\in \CC^{p}$ and an exact sequence in $\CC^{p}$
\begin{align*}
0 \to H^{-1}(E)[1] \to E \to H^{0}(E) \to 0. 
\end{align*}
The sheaf $H^{-1}(E)$ fits into the short exact sequence of sheaves
\begin{align*}
0 \to T_1 \to H^{-1}(E) \to T_2 \to 0, 
\end{align*}
with $T_1 \in \TT_{\omega, B}$ and $T_2 \in \FF_{\omega, B}$.
By the above two sequences, we obtain the exact sequences in $\CC^{p}$
\begin{align}\label{TEU}
&0 \to T_1[1] \to E \to U \to 0, \\
\label{TUH}
&0 \to T_2[1] \to U \to H^{0}(E) \to 0, 
\end{align}
for some $U\in \CC^{p}$.
Since $T_1\in \BB_{\omega, B}$
and $U \in \BB_{\omega, B}$ by the sequence (\ref{TUH}), 
the sequence (\ref{TEU}) implies 
\begin{align*}
H_{\BB}^{-1}(E)=T_1 \qquad \text{ and } \qquad H_{\BB}^{0}(E)=U. 
\end{align*}
In order to conclude $E\in \AA_{m\omega, B}$ for $m\gg 0$, it is enough to show that
\begin{align*}
\nu_{m\omega, B; \rm{min}}(T_1) \le 0 \qquad \text{ and } \qquad \nu_{m\omega, B; \rm{min}}(U)>0,
\end{align*}
for $m\gg 0$.
We only show that $\nu_{m\omega, B; \rm{min}}(U)>0$.
The inequality $\nu_{m\omega, B; \rm{min}}(T_1) \le 0$ is similarly proved, and we omit the proof.

Let $U^{m} \in \BB_{\omega, B}$ be the $\nu_{m\omega, B}$-semistable factor of $U$ such that
\[
\nu_{m\omega, B; \rm{min}}(U)=\nu_{m\omega, B}(U^{m}).
\]
We have the exact sequence in $\BB_{\omega, B}$
\begin{align*}
0 \to U^{'m} \to U \to U^m \to 0,
\end{align*}
and the long exact sequence of coherent sheaves
\begin{equation}\label{seq:L1}
\begin{split}
0 &\to H^{-1}(U^{'m}) \to T_2 \to H^{-1}(U^m) \\
& \to H^{0}(U^{'m}) \to H^{0}(E) \to H^{0}(U^m) \to 0.
\end{split}
\end{equation}
Since $H^{0}(E) \in \TT^{p}$, the sheaf $H^{0}(U^{m})$ also satisfies $H^0(U^m) \in \TT^p$.
This implies the inequality
\begin{align}\label{ineq:L1}
\nu_{m\omega, B}(H^{0}(U^{m}))>0, \ m>0.
\end{align}
Next we see that $\nu_{m\omega, B}(H^{-1}(U^{m})[1])$ is positive for $m\gg 0$.
The sequence (\ref{seq:L1}) gives rise to two short exact sequences of coherent sheaves
\begin{align*}
& 0\to K \to H^{-1}(U^m) \to K' \to 0, \\
& 0 \to K' \to H^{0}(U^{'m}) \to K'' \to 0.
\end{align*}
Note that there is a surjection $T_2 \twoheadrightarrow K$, hence $\mu_{\omega, B; \rm{min}}(K)$ is
bounded below, i.e., there is a constant $c$ which does not depend on $m$ such that $\mu_{\omega, B; \rm{min}}(K) \ge c$.
Also since $K''$ is a subsheaf of $H^{0}(E)$, it is a torsion sheaf and its first Chern class is bounded above.
This fact, together with $\mu_{\omega, B; \rm{min}}(H^0(U^{'m}))>0$, easily implies that
$\mu_{\omega, B; \rm{min}}(K')$ is bounded below.
Therefore $\mu_{\omega, B; \rm{min}}(H^{-1}(U^m))$ is also bounded below.

Let $A_1, \cdots, A_N$ be the set of 
 $\mu_{\omega, B}$-semistable factors of $H^{-1}(U^m)$.
From what we have proved above, there is a constant $c>0$, which does not depend on $m$ such that
\begin{align*}
0 \ge \mu_{\omega, B}(A_1) > \cdots >\mu_{\omega, B}(A_N) \ge -c.
\end{align*}
By replacing $c$ if necessary, the Hodge Index Theorem and the Bogomolov-Gieseker inequality imply the following bound:
\begin{align*}
\omega \cdot \frac{\cht_2(A_i)}{\cht_0(A_i)}  \le \omega \cdot \frac{\cht_1(A_i)^2}{2 \cht_0(A_i)^2} \le c.
\end{align*}
Hence we have
\begin{align}\label{ineq:L2}
\omega \cdot \frac{\cht_2 H^{-1}(U_m)}{\cht_0 H^{-1}(U_m)} \le c. 
\end{align}
Now we have
\begin{align*}
\frac{\Im Z_{m\omega, B}(H^{-1}(U^m)[1])}{\cht_0(H^{-1}(U^m))}=\frac{1}{6}m^3 \omega^3 -m\omega\frac{\cht_2 H^{-1}(U^m)}{\cht_0 H^{-1}(U^m)},
\end{align*}
which is positive for $m\gg 0$ by \eqref{ineq:L2}.
This implies the inequality
\begin{align}\label{ineq:L3}
\nu_{m\omega, B}(H^{-1}(U^m)[1])>0, \ m \gg 0.  
\end{align}
By the inequalities \eqref{ineq:L1} and \eqref{ineq:L3}, we obtain $\nu_{m\omega, B}(U^m)>0$ for $m\gg 0$.
\end{Prf}

Let us denote $\sigma_m=(Z_{m\omega, B}, \AA_{m\omega, B})$.
We have the following proposition.

\begin{Prop} \label{prop:large-volume}
Suppose that $\sigma_m$ are stability conditions.
If an object $E\in \Db(X)$ is $\sigma_m$-semistable for $m \gg 0$, then $E$ is semistable with
respect to the polynomial stability condition $(Z_{\infty\omega, B}, \CC^p)$.
\end{Prop}

\begin{Prf}
We may assume that $E\in \AA_{m\omega, B}$ for $m\gg 0$.
By Lemma~\ref{Lem:L1}, we have $E\in \CC^{p}$.
Suppose that $E$ is not semistable w.r.t. $Z_{\infty\omega, B}$.
Then there is an exact sequence in $\CC^p$
\begin{align}\label{seq:L2}
0 \to E' \to E \to E'' \to 0, 
\end{align}
such that $\arg Z_{m\omega, B}(E')> \arg Z_{m\omega, B}(E'')$ for $m\gg 0$.
However, by Lemma~\ref{Lem:L2}, the sequence \eqref{seq:L2} is also an exact sequence in
$\AA_{m\omega, B}$, for $m\gg 0$.
This implies that \eqref{seq:L2} destabilizes $E$ w.r.t. $\sigma_m$ for $m\gg 0$, which is a contradiction.
\end{Prf}

\section{Bogomolov-Gieseker type inequalities}\label{sec:inequality}

In this section we discuss bounds on the set of numerical classes of $\nu_{\omega,B}$-semistable objects in $\BB_{\omega, B}$.

\subsection{Torsion sheaves}\label{subsec:TorsionSheaves}
We consider Conjecture~\ref{Con:stability} in the case of torsion sheaves.
For simplicity, we assume that
\begin{align*}
\mathrm{Pic}(X)=\mathbb{Z}[\OO_X(H)],
\end{align*}
for an ample divisor $H$.
We denote $d=H^3 \in \mathbb{Z}_{>0}$ and $\omega=\alpha H$ for $\alpha \in \mathbb{Q}_{>0}$.
We take a smooth divisor
\begin{align*}
S \in \lvert mH \rvert,
\end{align*}
for $m \in \mathbb{Z}_{\ge 1}$ and consider semistable sheaves on $S$.
Notice that 
\begin{align*}
\Coh S \subset \BB_{\omega,B}, 
\end{align*}
is an abelian subcategory, and
for $E \in \Coh S$, we have $\nu_{\omega, B}(E)=\widehat\mu_{\omega, B}(E)$,
where $\widehat\mu_{\omega, B}$ is as defined in equation \eqref{eq:def-muhat}.
Hence if $E \in \Coh S$ is $\nu_{\omega,B}$-(semi)stable in $\BB_{\omega,B}$, 
then it is also
$\widehat{\mu}_{\omega,B}$-(semi)stable in $\Coh S$.
Let us discuss Conjecture~\ref{Con:stability} for objects in $\Coh S$.

Let $E\in \Coh(S)$ be a $\widehat{\mu}_{\omega,B}$-semistable sheaf with
\begin{align*}
\cht(E)=(r, l, s) \in H^0(S) \oplus H^2(S) \oplus H^4(S),
\end{align*}
where $\cht$ on $S$ is nothing but $\ch\cdot e^{-B|_S}$.
Let $i \colon S \hookrightarrow X$ denote the inclusion.
By the Grothendieck Riemann-Roch formula, we have
\begin{align*}
\cht(i_{\ast}E) &=\left(0, rS, -\frac{r}{2}S^2 +i_{\ast}l, \frac{r}{6}S^3 -\frac{1}{2}S\cdot i_{\ast}l +s \right)  \\
& \quad \in H^0(X) \oplus H^2(X) \oplus H^4(X) \oplus H^6(X).
\end{align*}
The condition $\nu_{\omega,B}(i_{\ast}E)=0$ is equivalent to $i_{\ast}l=rS^2/2$; hence we have
\begin{align*}
\cht(i_{\ast}E)= \left( 0, rS, 0, s-\frac{1}{12}rS^3 \right).
\end{align*}
On the other hand, the Bogomolov-Gieseker inequality $l^2 \ge 2rs$ implies the inequality
\begin{align}\label{Bog1}
s \le \frac{1}{8}rS^3 =\frac{1}{8}dr m^3.
\end{align}
Therefore we have
\begin{align*}
\Re Z_{\omega,B}(i_{\ast}E)
&= -s + \frac{1}{12}rS^3 +\frac{1}{2}\alpha^2 r H^2 S \\
&\ge \frac{rdm}{24} \left( 12\alpha^2 -m^2  \right).
\end{align*}
If $m>2\sqrt{3}\alpha$, we cannot conclude that $\Re Z_{\omega,B}(i_{\ast}E)$ is positive.
This implies that $\hat{\mu}_{\omega,B}$-stability and the Bogomolov-Gieseker inequality on $S$ are not sufficient to conclude Conjecture~\ref{Con:stability}, and we need to investigate $\nu_{\omega,B}$-stability in more detail.

Suppose, for instance, that there exists a $\mu_{\omega,B}$-semistable torsion free sheaf $F$ on $X$ such that $F|_{S} \cong E$.
For example, when $E=\OO_S(mH/2)$ for an even number $m$, we can take $F=\OO_X(mH/2)$.
We have
\begin{align*}
\mu_{\omega,B}(F)=\frac{\alpha^2}{2}md>0, \quad \mu_{\omega,B}(F(-S))=-\frac{\alpha^2}{2}md<0.
\end{align*}
Then the sequence
\begin{align*}
F \to E \to F(-S)[1]
\end{align*}
is an exact sequence in $\BB_{\omega,B}$.
If $E$ is $\nu_{\omega,B}$-semistable, we have $\nu_{\omega,B}(F) \le 0$.
Therefore we have
\begin{align}\label{Bog2}
s \le \frac{1}{6}rdm\alpha^2.
\end{align}
Note that, when $m$ is big, \eqref{Bog2} is a stronger inequality than \eqref{Bog1}.
Using (\ref{Bog2}) instead of (\ref{Bog1}), we obtain
\begin{align*}
\Re Z_{\omega,B}(i_{\ast}E) \ge \frac{rdm}{12} (4\alpha^2 +m^2)>0.
\end{align*}

By considering at the same time the two inequalities \eqref{Bog1} and \eqref{Bog2}, we obtain 
the inequality of Conjecture \ref{con:strong-BG} in this case:

\begin{Prop}\label{prop:SupportedSheaves}
Let $X$ be a smooth projective threefold with $\mathrm{Pic}(X)=\mathbb{Z}[\OO_X(H)]$.
Let $S\in\lvert mH \rvert$ be a smooth divisor, $i\colon S\into X$.
Let $E\in\Coh(S)$ be such that:
\begin{itemize}
\item $i_*E$ is $\nu_{\omega,B}$-semistable with $\nu_{\omega,B}(i_*E)=0$, and
\item there exists a torsion-free $\mu_{\omega,B}$-semistable sheaf $F\in\Coh(X)$ with $F|_S\cong E$.
\end{itemize}
Then
\begin{equation}\label{eq:Bonn06022011}
\cht_3(i_*E)\leq\frac{\omega^2}{18}\cht_1(i_*E).
\end{equation}
\end{Prop}
\begin{Prf}
First of all, when $3m^2\leq4\alpha^2$, we use \eqref{Bog1}, and we have
\[
\cht_3(i_*E)=s-\frac 1{12} rm^3d\leq\frac{rm^3d}{24}\leq\frac{r\alpha^2md}{18}=\frac{\omega^2}{18}\cht_1(i_*E).
\]
Similarly, when $3m^2\geq4\alpha^2$, we use \eqref{Bog2}, and \eqref{eq:Bonn06022011} is proved.
\end{Prf}

\subsection{Semistable objects in $\BB_{\omega, B}$}\label{subsec:SemistableObjectsInB}
Here we fix $B\in \NS_{\mathbb{Q}}(X)$ and use the twisted Chern character $\cht(E) =\ch (E) \cdot e^{-B}$.
For an ample divisor $\omega \in \NS_{\mathbb{Q}}(X)$, let
$\DDD \subset \BB_{\omega, B}$ denote the set
of objects $E \in \BB_{\omega, B}$ satisfying one of the following conditions:
\begin{enumerate}
\item $H^{-1}(E)=0$ and $H^0(E)$ is a pure sheaf of dimension $\ge 2$ that is
slope-semistable with respect to $\omega$.
\item $H^{-1}(E)=0$ and $H^0(E)$ is a sheaf of dimension $\le 1$.
\item \label{enum:D-Hm1}
$H^{-1}(E)$ is a torsion-free slope-semistable sheaf and $H^0(E) \in \Coh^{\le 1} X$.
If $\mu_{\omega, B}(H^{-1}(E))<0$, we have
\begin{equation}\label{eqn:S1}
\Hom(\Coh^{\le 1} X, E)=0.
\end{equation}
\end{enumerate}

\begin{Lem}\label{lem:large}
If an object $E\in \BB_{\omega, B}$ is $\nu_{m\omega, B}$-semistable for $m\gg 0$ then $E\in \DDD$.
\end{Lem}
\begin{Prf}
The proof is given by the same argument as in~\cite[Lemma~4.2]{large-volume}.
\end{Prf}

The set of objects $\DDD$ is also obtained as $\nu_{\omega, B}$-semistable objects $E$ with small $\omega^2 \cht_1(E)$.
Because of the rationality of $B$ and $\omega$, we can take $c\in \mathbb{Q}_{>0}$ to be
\begin{align}\label{minch}
c:= \mathrm{min} \left\{ \omega^2 \cht_1(E)>0 : E \in \BB_{\omega, B} \right\}.
\end{align}

\begin{Lem}\label{c1=H}
For an object $E\in \BB_{\omega, B}$, suppose that $\omega^2 \cht_1(E) \le c$.
Then $E$ is $\nu_{\omega, B}$-semistable if and only if $E\in \DDD$.
\end{Lem}

\begin{Prf}
The assumption implies that $\omega^2 \cht_1(E)=0$ or $c$, and the first case is obvious.
So assume $\omega^2 \cht_1(E)=c$. Given any short exact sequence
\begin{equation} \label{eq:AEB}
A \into E \onto B
\end{equation}
in $\BB_{\omega, B}$, we always have either $\omega^2 \cht_1(A) = 0$
or $\omega^2 \cht_1(B) = 0$. In the former case, we have
$\nu_{\omega, B}(A) = +\infty > \nu_{\omega, B}(E) \neq \infty$ and thus $E$ is unstable.
In the latter case, we have
$\nu_{\omega, B}(E) = +\infty > \nu_{\omega, B}(E)$, and thus the short
exact sequence \eqref{eq:AEB} cannot destabilize $E$.

So $E$ is stable if and only if there is no
subobject $A \into E$ with $\omega^2 \cht_1(A) = 0$, i.e $A$ is contained in the subcategory
described in \eqref{026}. This condition
is invariant under replacing $\omega$ by a scalar multiple $m \omega $ for $m \in \R$.
In particular, if $E$ is $\nu_{\omega, B}$-semistable, then it is
$\nu_{m\omega, B}$-semistable for $m \gg 0$, and $E \in \DDD$ by
Lemma \ref{lem:large}. 

Conversely, assume $E \in \DDD$ has a subobject $A$
with $\omega^2 \cht_1(A) = 0$. If $H^{-1}(A) \neq 0$, then 
it has slope $\mu_{\omega, B}(H^{-1}(A)) = 0$. Also, $H^{-1}(E)$ is non-trivial,
and thus $E$ satisfies condition \eqref{enum:D-Hm1} in the definition of $\DDD$.
As $\omega^2\cht_1(E) > 0$, we have $\mu_{\omega, B}(H^{-1}(E)) < 0$; but
$H^{-1}(A)$ is a subobject of $H^{-1}(E)$, in contradiction to the slope-semistability 
of $E$.

On the other hand, if $H^{-1}(A) = 0$, then $A \in \Coh^{\le 1} X$, in contradiction to \eqref{eqn:S1}.
\end{Prf}

In particular, Conjecture~\ref{con:strong-BG} includes the following conjecture for
$\mu_{\omega, B}$-stable sheaves:

\begin{Con}\label{conj:ineq:sheaf}
Let $X$ be a smooth projective threefold and take $B, \omega \in \NS_{\mathbb{Q}}(X)$ with $\omega$ ample.
Let $E$ be a $\mu_{\omega, B}$-stable sheaf satisfying $\omega^2 \cht_1(E)=c$, where $c$ is defined by \eqref{minch}.
Suppose that $E$ satisfies
\begin{align}\label{eq:sheaf1}
\frac{\omega \cht_2(E)}{\omega^3 \cht_0(E)}=\frac{1}{6}.
\end{align}
Then we have 
\begin{align}\label{ineq:sheaf2}
\frac{\cht_3(E)}{\omega^2 \cht_1(E)} \le \frac{1}{18}.
\end{align}
\end{Con}

\begin{Ex} \label{ex:L-otimes-IC}
In the situation of Conjecture~\ref{conj:ineq:sheaf}, 
suppose that $\mathrm{Pic}(X)$ is generated by 
an ample line bundle $L$ on $X$
with $D=L^3$. 
Setting $\omega =tL$ for $t\in \mathbb{Q}_{>0}$
and $B=0$,  
we have $c=Dt^2$. 
For a curve $C \subset X$ of degree $d = L.C$, let
$I_C$ be the ideal sheaf of $C$.
Then
the object $E=L\otimes I_C$ 
is $\mu_{\omega, 0}$-stable with $\omega^2 \cht_1(E)=c$, 
and satisfies (\ref{eq:sheaf1}) if $d < \frac D2$ and
\begin{align*}
\frac{t^2}6 = \frac 12 - \frac dD.
\end{align*}
Then the Conjecture states that 
\begin{align}
\ch_3(E) =  \frac D6 - d - \ch_3(\OO_C) & \le \frac{t^2D}{18} = \frac D6 - \frac d3  \nonumber\\
\intertext{or, equivalently,}
-\ch_3(\OO_C) & \le \frac 23 d. \label{ineq:ch3}
\end{align}
For instance, if $X \subset \P^4$ is a hypersurface of degree $D$, Hirzebruch-Riemann-Roch
relates $\ch_3(\OO_C)$ to the arithmetic genus $g$ of $C$ by
\begin{align*}
1 - g = \chi(\OO_C) = \ch_3(\OO_C) + \frac d2 (4 - D).
\end{align*}
Thus the inequality \eqref{ineq:ch3} becomes
$ g \le \frac d2 D - \frac 43 d + 1$.
Since $d < \frac D2$, this follows from Castelnuovo's classical inequality
$ g \le \frac 12 (d-1)(d-2)$, which does hold in our situation: it has been shown for
for singular curves $C \subset \P^3$ in \cite{Okonek-Spindler:spektrum-II}
and \cite{Hartshorne:space-curves}; and since $C \subset X$ is contained in a smooth hypersurface
in $\P^4$, the curve $C$ will map isomorphically into $\P^3$ under a generic projection.

On the other hand, already when $X \subset \P^{N+3}$ is a complete intersection of codimension
$N$, the inequality \eqref{ineq:ch3} seems stronger than known Castelnuovo
inequalities: it becomes
\[
g \le \frac{D_1 + \dots + D_N}2 d - \frac{N+3}2 d + \frac 23 d + 1
\]
for any curve of genus $g$ and degree $d < \frac 12 D_1 \cdot D_2 \cdots D_N$ on a complete
intersection of degree $(D_1, \dots, D_N)$. The statement would be similar to the case of space
curves lying on a surface of given degree, treated e.g.~ in \cite{Harris:space-curves}.
\end{Ex}

\subsection{Bogomolov-Gieseker type inequality without $\ch_3$} \label{sec:BGwithoutch3}
In this section we establish a Bogomolov-Gieseker type inequality for $\nu_{\omega, B}$-semistable
objects in $\BB_{\omega, B}$ which does not involve $\ch_3$.
For $a, b\in \mathbb{R}$, we set $f_{a, b} \colon \NS_{\mathbb{Q}}(X) \to \mathbb{R}$ to be
\begin{align*}
f_{a, b}(x) := a\omega^3 \cdot x^2 \omega + b(x\omega^2)^2. 
\end{align*}
Recall that the discriminant $\Delta(E) \in A^2_{\mathbb{Q}}(X)$ of an object $E\in \Db(X)$ is defined by
\begin{align*}
\Delta(E) &=\ch_1(E)^2 -2\ch_0(E) \ch_2(E) \\
&=\cht_1(E)^2 -2\cht_0(E)\cht_2(E).
\end{align*}

We have the following result.

\begin{Thm}\label{Thm:Bog}
For a smooth projective threefold $X$ and $\omega \in \NS_{\mathbb{Q}}(X)$, take $a \in \mathbb{R}_{\ge -1}$ and $b\in \mathbb{R}$ such that $f_{a, b}$ satisfies the following conditions:
\begin{enumerate}
\item \label{item1}$f_{a, b}(x) \ge 0$, for any $x\in \NS_{\mathbb{Q}}(X)$.
\item \label{item2} $f_{a+1, b}(x) \ge 0$, for any effective class $x\in \NS_{\mathbb{Q}}(X)$.
\end{enumerate}
Then, for any $\nu_{\omega, B}$-semistable object $E\in \BB_{\omega, B}$, we have the following inequality:
\begin{align}\label{ineq:Bog}
\omega^3 \cdot \omega \Delta(E) +f_{a, b}(\cht_1(E)) \ge 0.
\end{align}
\end{Thm}

\begin{Prf}
We prove the inequality (\ref{ineq:Bog}) by induction on $\omega^2\cht_1(E)$.
Observe that the Bogomolov-Gieseker inequality and the condition \eqref{item1} imply
\eqref{ineq:Bog} for torsion free slope-semistable sheaves.
Also, condition~\eqref{item2} implies the inequality \eqref{ineq:Bog} for torsion sheaves.
Therefore \eqref{ineq:Bog} holds for any object $E\in \DDD$; hence it also holds if $\omega^2 \cht_1(E) \le c$, by Lemma~\ref{c1=H}.

Assume that (\ref{ineq:Bog}) holds for all $\nu_{\omega, B}$-semistable $F \in \BB_{\omega, B}$ with $\omega^2 \cht_1(F)< \omega^2 \cht_1(E)$.
By the previous argument, we may assume that $E\notin \DDD$.
Then, by Lemma~\ref{lem:large}, the object $E$ is not $\nu_{m\omega, B}$-semistable for sufficiently large $m$.
Hence we can take $m_0 \in \mathbb{R}_{>0}$ to be
\begin{align*}
m_0=\mathrm{sup} \{ m\in \mathbb{R}_{>0} \colon E \mbox{ is }\nu_{m\omega, B}\mbox{-semistable } \}.
\end{align*}
By Corollary \ref{cor:tilt-wall-crossing}, there is a filtration in $\BB_{\omega, B}$
\begin{align*}
0=E_0 \subset E_1 \subset \cdots \subset E_N=E, 
\end{align*}
such that the following holds.
\begin{itemize}
\item $E_{\bullet}$ is a Harder-Narasimhan filtration of $E$ with respect to
$\nu_{(m_0+\varepsilon)\omega, B}$-stability for $0<\varepsilon \ll 1$.
In particular, subquotients $F_i=E_i/E_{i-1}$ satisfy
\begin{align}\label{ineq:nu1}
\nu_{(m_0+\varepsilon)\omega, B}(F_1)>\nu_{(m_0+\varepsilon)\omega, B}(F_2)>\cdots
>\nu_{(m_0+\varepsilon)\omega, B}(F_N).
\end{align}
\item The subquotients $F_i$ are $\nu_{m_0 \omega, B}$-semistable with
\begin{align}\label{eq:nu1}
\nu_{m_0\omega, B}(F_1)=\nu_{m_0\omega, B}(F_2)=\cdots =\nu_{m_0\omega, B}(F_N).
\end{align}
\end{itemize}
We set $a_i$, $b_i$ and $c_i$ as follows:
\begin{align*}
a_i =\frac{\omega^3 \cht_0 (F_i)}{\omega^2 \cht_1(F_i)}, \ b_i =\frac{\omega \cht_2 (F_i)}{\omega^2 \cht_1(F_i)}, \ c_i=\frac{\cht_1 (F_i)}{\omega^2 \cht_1(F_i)}.
\end{align*}
Note that we have $\omega^2 \cht_1(F_i) < \omega^2 \cht_1(E)$; hence $F_i$ satisfies the inequality \eqref{ineq:Bog} by the inductive assumption.
The inequality \eqref{ineq:Bog} for $F_i$ is written as
\begin{align}\label{ineq:induction}
(a+1)\omega^3 \cdot \omega c_i^2-2a_i b_i +b \ge 0.
\end{align}
By setting $c=m_0^2/6>0$, the equality \eqref{eq:nu1} implies
\begin{align}\label{eq:nu2}
-c a_1+b_1=-c a_2 +b_2= \cdots =-c a_N +b_N.
\end{align}
Combined with \eqref{ineq:nu1}, we obtain the inequalities
\begin{align}\label{ineq:nu2}
a_1<a_2< \cdots <a_N.
\end{align}
Then \eqref{eq:nu2} and \eqref{ineq:nu2} imply the following inequalities:
\begin{align}\label{ineq:nu3}
b_1<b_2< \cdots <b_N.
\end{align}
We can calculate as
\begin{align*}
&\omega^3 \cdot \omega \Delta(E) +f_{a, b}(\cht_1(E)) \\
&=\omega^3 \cdot \omega\left( \left( \sum_{i=1}^{N} \cht_1(F_i) \right)^2 -2 \left(\sum_{i=1}^{N} \cht_0(F_i) \right) \left( \sum_{i=1}^{N} \cht_2(F_i)\right)  \right) \\
& \hspace{10mm} +a\omega^3 \cdot \omega \left( \sum_{i=1}^{N} \cht_1(F_i) \right)^2 + b \left(\sum_{i=1}^{N} \omega^2 \cht_1(F_i) \right)^2 \\
&= \sum_{i=1}^{N} \omega^3 \cdot \omega \Delta(F_i)+f_{a, b}(\cht_1(F_i)) \\
& \hspace{5mm} +2\sum_{i<j} \omega^2 \cht_1(F_i) \cdot \omega^2 \cht_1(F_j) \left((a+1)\omega^3 \cdot \omega c_i c_j - a_i b_j-a_j b_i +b  \right).
\end{align*}
The first term of the last equation is non-negative by the inductive assumption.
As for the second term, note that (\ref{ineq:induction}) implies
\begin{align*}
b \ge -\frac{1}{2}(a+1)\omega^3 \cdot (\omega c_i^2 +\omega c_j^2)+a_ib_i +a_j b_j.
\end{align*}
Therefore we have
\begin{align}\notag
&(a+1)\omega^3 \cdot \omega c_i c_j - a_i b_j-a_j b_i +b \\
\label{ineq:last} &\ge (a_j-a_i)(b_j-b_i) -\frac{1}{2}(a+1)\omega^3 \cdot \omega(c_i-c_j)^2.
\end{align}
Note that since $\omega^2(c_i-c_j)=0$, the Hodge Index Theorem implies $\omega(c_i-c_j)^2 \le 0$.
Combined with $a\ge -1$, \eqref{ineq:nu2} and \eqref{ineq:nu3}, we conclude that $\eqref{ineq:last} \ge 0$.
By induction, we obtain the desired inequality \eqref{ineq:Bog}.
\end{Prf}

\begin{Cor}\label{Cor:Bog1}
Let $X$ be a smooth projective threefold and take $B, \omega \in \NS_{\mathbb{Q}}(X)$ with $\omega$ ample.
Then, for any $\nu_{\omega, B}$-semistable object $E\in \BB_{\omega, B}$, we have
\begin{align*}
\overline{\Delta}_{\omega}(E):=(\omega^2 \cht_1(E))^2 -2\omega^3 \cht_0(E) \cdot \omega\cht_2(E) \ge 0.
\end{align*}
\end{Cor}
\begin{Prf}
We take $a=-1$ and $b=1$ in Theorem~\ref{Thm:Bog}.
The condition (\ref{item1}) is satisfied by the Hodge Index Theorem, and \eqref{item2} is obvious.
Therefore the result follows from Theorem~\ref{Thm:Bog}.
\end{Prf}

\begin{Cor}\label{Cor:Bog2}
Under the previous assumptions, there is a constant $C_{\omega} \in \mathbb{R}_{\ge 0}$, which depends
only $[\omega] \in \mathbb{P}(\NS_{\mathbb{Q}}(X))$, such that any $\nu_{\omega, B}$-semistable
object $E\in \BB_{\omega, B}$ satisfies the following inequality:
\begin{align*}
\omega^3 \cdot \omega \Delta(E) +C_{\omega}(\omega^2 \cht_1(E))^2 \ge 0.
\end{align*}
\end{Cor}
\begin{Prf}
We set $a=0$ and want to find $b=C_{\omega} \ge 0$ such that the conditions \eqref{item1} and \eqref{item2} are satisfied.
The condition \eqref{item1} is obviously satisfied for any $b\ge 0$, so it is enough to deal with
\eqref{item2}: it requires that for any effective divisor $D$, we have
\begin{align}\label{ineq:D}
\omega^3 \cdot D^2 \omega +C_{\omega} (D\omega^2)^2 \ge 0
\end{align}
holds for any effective divisor $D$.

Fix a norm $\lVert \ast \rVert$ on $H^2(X, \R)$. Then there is a constant
$A_\omega$ such that
\[
\omega^3 \cdot \omega D^2 \le A_\omega \lVert D \rVert^2
\]
for every $D \in H^2(X, \R)$.
On the other hand, due to the openness of the ample cone, there is a constant
$B_\omega$ such that
\[
\omega^2 D \ge B_\omega \lVert D \rVert
\]
for every effective divisor class $D$. Setting $C_\omega := \frac{A_\omega}{B_\omega^2}$, we obtain
\eqref{ineq:D} as required.
\end{Prf}

\begin{Cor}\label{Cor:Bog3}
Let $X$ be a smooth projective threefold such that $\NS(X)$ has rank one.
Then for any $\nu_{\omega, B}$-semistable object $E\in \BB_{\omega, B}$, we have
\begin{align*}
\omega \Delta(E) \ge 0.
\end{align*}
\end{Cor}
\begin{Prf}
If $\NS(X)$ is of rank one, then we can take $a=b=0$ in Theorem~\ref{Thm:Bog}.
\end{Prf}

\subsection{Stability of vector bundles with trivial discriminant}\label{subsec:LineBundles}
As an application of Corollary~\ref{Cor:Bog1}, we have the following result which generalizes \cite[Prop.\ 3.6]{AB:Reider}.

\begin{Prop}\label{prop:StabilityLineBundles}
Let $E$ be a $\mu_{\omega, B}$-stable vector bundle on $X$ with $\overline{\Delta}_{\omega}(E)=0$.
Then $E$ is $\nu_{\omega, B}$-stable.
\end{Prop}
\begin{Prf}
By Proposition \ref{prop:ZB-selfdual} and Section \ref{sec:Comparison}, we may replace
$E$ by its dual, and thus may assume $\omega^2\cht_1(E)\leq0$; in particular, $E[1]\in\BB_{\omega, B}$.
Assume, for a contradiction, $E[1]$ is not $\nu_{\omega, B}$-stable.
Then $\omega^2\cht_1(E)<0$, and there exists a destabilizing sequence
\[
M\to E[1]\to N[1],
\]
with $N$ a $\nu_{\omega, B}$-semistable sheaf such that
\begin{equation}\label{eq:starting}
\nu_{\omega, B}(N[1])\leq\nu_{\omega, B}(E[1]).
\end{equation}
Expanding \eqref{eq:starting}, by using the assumption that $\overline{\Delta}_{\omega}(E)=0$, we deduce the inequality
\[
\frac{2\omega\cht_2(N)}{\omega^2\cht_1(N)}-\frac{\omega^3\rk(N)}{3\omega^2\cht_1(N)}\leq\frac{\omega^2\cht_1(E)}{\omega^3\rk(E)}-\frac{\omega^3\rk(E)}{3\omega^2\cht_1(E)}.
\]
Now, $\omega^2\cht_1(N)<0$ and Corollary \ref{Cor:Bog1} give
\[
\frac{\omega^2\cht_1(N)}{\omega^3\rk(N)}-\frac{\omega^3\rk(N)}{3\omega^2\cht_1(N)}\leq\frac{\omega^2\cht_1(E)}{\omega^3\rk(E)}-\frac{\omega^3\rk(E)}{3\omega^2\cht_1(E)}.
\]
Therefore, we get the inequality
\[
\frac{1}{\omega^3}\left( \mu_{\omega,B}(E)-\mu_{\omega,B}(N)\right)\geq\frac{\omega^3}{3}\left(\frac{\mu_{\omega,B}(N)-\mu_{\omega,B}(E)}{\mu_{\omega,B}(N)\mu_{\omega,B}(E)} \right).
\]
But $\mu_{\omega,B}(E)<\mu_{\omega,B}(N)$: Indeed, otherwise, the $\mu_{\omega,B}$-stability of $E$ would imply that
$\Hom(H^0(M),Q)\neq0$, where $Q$ is a $\mu_{\omega, B}$-stable quotient of $N$ with slope $\leq\mu_{\omega,B}(E)$, violating Lemma \ref{lem:Cohomologies1}.
This gives a contradiction, since $\omega^3$ and $\mu_{\omega,B}(N)\mu_{\omega,B}(E)$ are both positive.
\end{Prf}

In particular, by Proposition \ref{prop:StabilityLineBundles}, all (shifts of) line bundles with $\overline{\Delta}_{\omega}=0$ are $\nu_{\omega, B}$-stable in $\BB_{\omega, B}$.
When $\NS(X)$ is of rank one, this is true for all line bundles on $X$.

Proposition \ref{prop:StabilityLineBundles} gives also a further evidence for Conjecture \ref{con:strong-BG} (and so for Conjecture \ref{Con:stability}).
More precisely, following Dr\'ezet (see e.g.~ \cite[Sect.\ 3.5]{Langer:Survey}), we introduce the higher discriminants $\overline{\Delta}_{\omega,i}$ as follows:
\begin{align*}
\overline{\Delta}_{\omega,1}&=\omega^2\cht_1\\
\overline{\Delta}_{\omega,2}&=\overline{\Delta}_{\omega}\\
\overline{\Delta}_{\omega,3}&=2\left(3(\omega^3\rk)^2\cht_3-3(\omega^3\rk)(\omega^2\cht_1)(\omega\cht_2)+(\omega^2\cht_1)^3\right).
\end{align*}
Notice that these higher discriminants are invariant under tensoring by line bundles whose numerical class is a multiple of $\omega$.

\begin{Prop}\label{prop:Simpson}
Let $E$ be a $\mu_{\omega, B}$-stable vector bundle on $X$ with $\overline{\Delta}_{\omega,2}(E)=0$.
Then $\overline{\Delta}_{\omega,3}(E)=0$.
\end{Prop}

\begin{Prf}
First of all, since $E$ is $\mu_{\omega,B}$-stable, we have $\omega\Delta(E)=0$ and so $(\omega^2\cht_1(E))^2=\omega^3(\omega\cht_1(E)^2)$.
Hence, by taking a finite cover and a tensor by a line bundle, we can reduce to the case $\omega^2\cht_1(E)=0$.
Our assumptions then give $\omega\cht_2(E)=0$ and
\[
\overline{\Delta}_{\omega,3}(E)=6(\omega^3\rk(E))^2\cht_3(E).
\]
We want to prove that $\cht_3(E)=0$.
Again, by taking a finite cover, we can assume that $B$ is the numerical class of a line bundle $L$.
Consider the $\mu_{\omega}$-stable vector bundle $F=E\otimes L^\vee$.
Then $\ch(F)=\cht(E)$.
In particular, we have
\[
\omega^2\ch_1(F) = \omega\ch_2(F)=0.
\]
By \cite[Thm.\ 2]{Simpson:Higgs} (see also \cite[Thm.\ 4.1]{Langer:S-group} for an algebraic proof), $\ch_3(F)=0$, and so $\cht_3(E)=0$, as wanted.
\end{Prf}

Let $E$ be as in the Proposition, and assume $\Im Z_{\omega,B}(E)=0$ and $\omega^2\cht_1(E)>0$.
Then Conjecture \ref{Con:stability2} is satisfied and the inequality of Conjecture \ref{con:strong-BG} becomes an equality:
\begin{equation*}
\begin{split}
\cht_3(E)-\frac{\omega^2}{2}\cht_1(E)&<\cht_3(E)-\frac{\omega^2}{18}\cht_1(E)\\
   &=\frac{1}{6(\omega^3\rk(E))^2}\left(\overline{\Delta}_{\omega,3}(E)-2
\overline{\Delta}_{\omega,1}(E)\overline{\Delta}_{\omega,2}(E)\right)=0.
\end{split}
\end{equation*}

\section{Examples}\label{sec:Examples}

In this section we discuss Conjecture \ref{Con:stability} in some examples, focusing on the case of the projective space.

\subsection{A technical result}\label{subsec:tech}
Let $X$ be a smooth projective threefold and let $B,\omega\in \NS_\Q(X)$ with $\omega$ ample.
We consider slightly more central central charges on $\AA_{\omega,B}$: for $s\in\Q$, define
\begin{align*}
Z_{\omega,B,s}(-):= \left(-\cht_3(-)+s\omega^2\cht_1(-) \right)
+i\left(\omega \cht_2(-)-\frac{\omega^3}{6}\cht_0(-)  \right).
\end{align*}
We have $Z_{\omega,B}=Z_{\omega,B,1/2}$.

\begin{Prop}\label{prop:tech}
Let $\CC\subset\Db(X)$ be a heart of a bounded t-structure with the following properties:
\begin{enumerate}
\item\label{tech:1} there exist $\phi_0\in(0,1)$ and $s_0\in\Q$ such that
\[
Z_{\omega,B,s_0}(\CC)\subset \left\{r\exp(\pi\phi i)\,:\, r\geq0,\, \phi_0\leq\phi\leq\phi_0+1\right\}.
\]
\item\label{tech:2} $\CC\subset\langle\AA_{\omega,B},\AA_{\omega,B}[1]\rangle$.
\item\label{tech:3} for all $x\in X$, we have $k(x)\in\CC$ and, for all proper subobjects $C\into k(x)$ in $\CC$, $\Im Z_{\omega,B}(C)>0$.
\end{enumerate}
Then the pair $(Z_{\omega,B,s},\AA_{\omega,B})$ is a stability condition on $\Db(X)$, for all $s>s_0$.
\end{Prop}



\begin{Prf}
To simplify the notation, we put $Z_s=Z_{\omega,B,s}$ and $\AA=\AA_{\omega,B}$.
By Corollary \ref{cor:HNfiltrations} and Lemma \ref{lem:Access}, to prove that $(Z_{s},\AA)$ is a stability condition it will be enough to prove that $\Re Z_{s_0}(\TT^1) < 0$, where $\TT^1$ is the abelian subcategory of $\AA$, defined in Lemma \ref{lem:Access}, of objects in $\AA$ with $\Im Z_{\omega,B}=0$.
Assume, for a contradiction, this is not the case.
Then, by Lemma \ref{lem:Access}, there exists a simple object $F\in\TT^1$ with $\Re Z_{s_0}(F)\geq0$.
If $\Re Z_{s_0}(F)=0$, then, from $\omega^2\cht_1(F)<0$, we deduce that $\Re Z_{s}(F)<0$, for all $s>s_0$.
Hence, we can assume that there exists $F\in\TT^1$ with $\Re Z_{s_0}(F)>0$.

By  Proposition \ref{prop:tilt-properties}, assumption \eqref{tech:2} implies that $\CC$ is a tilt of $\AA$.
Consider the torsion pair on $\AA$ induced by $\CC$:
\begin{align*}
&\TT:=\AA\cap\CC\\
&\FF:=\AA\cap\left(\CC[-1]\right).
\end{align*}
By assumption \eqref{tech:1}, for all simple objects $E\in\TT^1$ with $Z_{s_0}(E)\neq0$, we have
\begin{itemize}
\item $E\in\TT$ if and only if $\Re Z_{s_0}(E)<0$.
\item $E\in\FF$ if and only if $\Re Z_{s_0}(E)>0$.
\end{itemize}
Indeed, by Lemma \ref{lem:Access}, a simple object in $\TT^1$ has no proper subobject in $\AA$, and so it belongs either to $\TT$ or to $\FF$.

Hence, $F\in\FF\cap\TT^1$.
Since $F$ is simple, $\Hom(F,k(x))=\Hom(k(x),F)=0$.
By Proposition \ref{prop:ZB-selfdual} and Remark \ref{rmk:cohomologies}, up to replacing $F$ with $F^\vee[3]$ (and $B$ with $-B$), we can assume that $\Hom(k(x),F[1])\neq0$, for all $x\in D$, where $D$ is a divisor in $X$.
Consider a non-zero morphism $f\colon k(x)\to F[1]$ and $\ker(f),\cok(f),\im(f)\in\CC$.
By assumption \eqref{tech:3}, if $\ker(f)\neq0$, then $\Im Z_{s_0}(\ker(f))>0$.
Hence, $\Im Z_{s_0}(\im(f))<0$, and so $\Im Z_{s_0}(\cok(f))>0$.
Since $\CC$ is a tilt of $\AA$, $\AA$ is a tilt of $\CC$ as well.
Consider the induced torsion pair on $\CC$:
\begin{align*}
&\TT':=\CC\cap\left(\AA[1]\right)\\
&\FF':=\CC\cap\AA.
\end{align*}
Then $F[1]\in\TT'$.
But $\TT'$ is closed under quotients.
This gives $\cok(f)\in\TT'$, and so $\Im Z_{s_0}(\cok(f))\geq0$, a contradiction.

Therefore, we have an exact sequence in $\CC$
\[
0\to k(x)\to F[1]\to Q_0[1]\to0,
\]
for some $Q_0[1]\in\CC$.
Since $Q_0$ is then an extension of $k(x)$ by $F$, it also belongs to $\AA$.
Hence, $Q_0\in\AA\cap(\CC[-1])=\FF$ and $Z_{s_0}(Q_0)=Z_{s_0}(F)-1$.
But, for $y\in D$, $y\neq x$, we have $\Hom(k(y),Q_0[1])\neq0$.
We can repeat the previous argument and construct a sequence of $Q_m\in\FF$, $m\in\N$, with
\[
Z_{s_0}(Q_m)=Z_{s_0}(Q_{m-1})-1=\ldots=Z_{s_0}(F)-m-1.
\]
But, for $m\gg0$, $Z_{s_0}(Q_m)<0$, a contradiction to $Q_m\in\FF$.
\end{Prf}

In the next section, we will apply Proposition \ref{prop:tech} together with Proposition \ref{prop:StabilityLineBundles} to give some examples in which Conjecture \ref{Con:stability} is verified.

\subsection{The projective space}\label{subsec:P3}
Consider the projective space $\P^3$.
For simplicity, let us fix $B = 0$, so that $\cht =\ch$.
Identifying $\Num_\Q(\P^3)$ with $\Q^{\oplus 4}$ in the obvious way, we define the central charge $Z^{s,t}$
for $s \in \Q, t \in \Q_{> 0}$ by
\[
Z^{s,t}(-)=\left(- \ch_3(-) + s \ch_1(-)\right) + i \left(\ch_2(-) - t \rk (-)\right).
\]
The central charge $Z_{\omega, 0}$ of equation \eqref{3fold-charge} corresponds to the choices
of $t=\frac{\omega^2}{6}$
$s=\frac{\omega^2}{2}$, up to an overall multiplication of the imaginary part of $Z_{\omega, B}$ by
$\omega^{-1}$. (The last operation is part of the $\widetilde{\GL}_2(\R)$-action on the set of
stability conditions defined in \cite{Bridgeland:Stab}, and does not affect the set of stable
objects.)

Given $t\in\Q$, consider the abelian category $\AA_t = \AA_{\omega, 0}$ for $t = \frac{\omega^2}6$
constructed as before---the only difference is that we only assume $\omega^2$ to be a rational number.
Conjecture \ref{Con:stability} then reads as follows: The pair $(Z^{s,t},\AA_t)$ is a Bridgeland
stability condition on $\Db(\P^3)$ for $s=3t$.
Our goal is to use Proposition \ref{prop:tech} to prove a strengthening of Conjecture \ref{Con:stability} in this case:

\begin{Thm}\label{thm:P3}
Let $s,t\in\Q$ be such that $0<t< \frac 12$ and
\begin{equation*}\label{eqn:Bonn4}
s>\frac{7 t-2}{6(t+1)}.
\end{equation*}
Then the pair $(Z^{s,t},\AA_t)$ defines a stability condition on $\Db(\P^3)$.
\end{Thm}

Since $s = 3t$ satisfies the above inequality for $0 < t < \frac 12$, this proves
Conjecture \ref{Con:stability} for $\omega<\sqrt{3}$. Moreover:

\begin{Rem}
The strong Conjecture \ref{con:strong-BG} for given $\omega$ and $B = 0$ holds if and only
if the pair $(Z^{s,t},\AA_t)$ defines a stability condition for
$t = \frac{\omega^2}6$ and $s = \frac t3$. Indeed, with this choice, a
tilt-stable objects $E \in \BB_{\omega, 0}$
with $\nu_{\omega, 0}(E) = 0$ satisfy
$Z^{s, t}(E) > 0$ if and only if $\ch_3(E) < s \ch_1(E) = \frac{\omega^2}{18} \ch_1(E)$.
Since $\frac t3 > \frac{7t-2}{6(t+1)}$ for $t < \frac 12$, Theorem \ref{thm:P3} actually proves
the strong Conjecture for $\omega < \sqrt{3}$.
\end{Rem}

Notice that, for $(s,t)=(1/6,1/2)$, we have $Z^{1/6,1/2}(\OO_{\P^3}(1))=0$, and so, by Lemma \ref{lem:algebraic} below, the function $Z^{1/6,1/2}$ does not define a stability condition.

To prove Theorem \ref{thm:P3}, recall that, by a classical result of Beilinson, on $\Db(\P^3)$ we have a bounded $t$-structure with heart given by
\[
\CC:=\langle \OO_{\P^3}(-1)[3], \OO_{\P^3}[2], \OO_{\P^3}(1)[1], \OO_{\P^3}(2)\rangle.
\]

An easy computation shows the following:

\begin{Lem}\label{lem:algebraic}
Assume $0<t<1/2$ and
\begin{equation}\label{eq:Bonn03022011}
\frac{7 t-2}{6(t+1)}<s\leq 1/6.
\end{equation}
Then the pair $\sigma_Q:=(Z^{s,t},\PP_Q)=(Z^{s,t},\CC)$ defines a stability condition on $\Db(\P^3)$.
The skyscraper sheaves $k(x)$, $x\in\P^3$, are $\sigma_Q$-stable of phase 1.
\end{Lem}

\begin{figure}[htb]
\begin{center}
\includegraphics{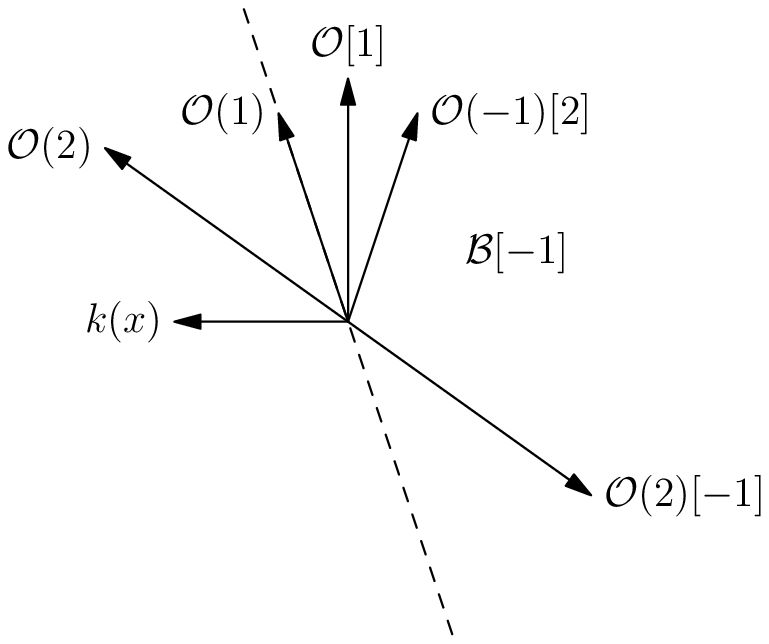}
\caption{Tilting $\AA$.}
\label{fig:tilt}
\end{center}
\end{figure}

\begin{Prf} (Theorem \ref{thm:P3})
Let $s_0$ be a rational number satisfying the inequalities \eqref{eq:Bonn03022011}.
We want to apply Proposition \ref{prop:tech} to $Z^{s_0,t}$, $\AA_t$, and $\CC$.
First of all, assumptions \eqref{tech:1} and \eqref{tech:3} in Proposition \ref{prop:tech} follow directly from Lemma \ref{lem:algebraic}, where $\phi_0:=\frac 1\pi\arg Z^{s,t}(\OO_{\P^3}(1))\in(0,1)$ is the phase of $\OO_{\P^3}(1)$.
Hence, we only need to show \eqref{tech:2}.
But, as an easy consequence of Proposition \ref{prop:StabilityLineBundles}, we have that the following objects are in $\AA_t$:
\begin{itemize}
\item $\OO_{\P^3} (k)$, for $k\geq1$,
\item $\OO_{\P^3}[1]$,
\item $\OO_{\P^3}(k)[2]$, for $k\leq-1$.
\end{itemize}
Hence, since $\CC$ is the category generated by extensions by  $\OO_{\P^3}(-1)[3]$, $\OO_{\P^3}[2]$, $\OO_{\P^3}(1)[1]$, and $\OO_{\P^3}(2)$, we have
\[
\CC\subset\langle \AA_t,\AA_t[1]\rangle,
\]
as wanted (see Figure \ref{fig:tilt}).
\end{Prf}

\begin{Rem}
Notice that, if $\widetilde{\AA}_t$ denotes the tilt of $\CC$ given by
\[
\widetilde{\AA}_t:=\PP_Q((0,1]),
\]
then, for $0<t<1/2$, we have $\AA_t=\widetilde{\AA}_t$.
Moreover, this shows the existence of stability conditions on $\Db(\P^3)$ for all irrational $t$ and $s$ satisfying inequalities \eqref{eq:Bonn03022011}.
\end{Rem}

\begin{Rem}\label{rmk:quadric}
The proof given in this section is in principle generalizable to other threefolds admitting a strong full exceptional collection.
An example is the quadric threefold $i_Q:Q\into\P^4$.
Denote by $S$ the spinor vector bundle of \cite{Kapranov:flagvarieties}, defined by an exact sequence
\[
0\to\OO_{\P^4}(-1)^{\oplus 4}\to\OO_{\P^4}^{\oplus 4}\to (i_Q)_*S\to 0.
\]
Then, for example, we have a strong full exceptional collection
\[
\left\{ \OO_Q(-1),S(-1),\OO_Q,\OO_Q(1)\right\},
\]
where $\OO_Q(1):=\OO_{\P^4}(1)|_Q$.
Identify $\mathrm{Pic}(Q)\cong\Z$ with generator the hyperplane section $h_Q$, and set $B=-\frac 12 h_Q$.
By considering the heart
\[
\CC:=\langle \OO_Q(-1)[3],S(-1)[2],\OO_Q[1],\OO_Q(1)\rangle,
\]
we have, by Proposition \ref{prop:tech}, that $Z_{\omega,B}$ is a stability condition for $\omega^3<\frac{1}{12\sqrt{3}}$.
(Indeed, all line bundles belong to $\AA_{\omega,B}$ up to shift by Proposition \ref{prop:StabilityLineBundles} and $S(-1)[1]\in\AA_{\omega,B}$ by our choice of $B$. The rest is precisely the same argument as in Theorem \ref{thm:P3}.)
\end{Rem}

\bibliography{all}                      

\begin{thebibliography}{BBMT11}

\bibitem[AB11]{AB:Reider}
Daniele Arcara and Aaron Bertram.
\newblock Reider's theorem and {T}haddeus pairs revisited.
\newblock In {\em Grassmannians, moduli spaces and vector bundles}, volume~14
  of {\em Clay Math. Proc.}, pages 51--68. Amer. Math. Soc., Providence, RI,
  2011.
\newblock arXiv:0904.3500.

\bibitem[ABL07]{Aaron-Daniele}
Daniele Arcara, Aaron Bertram, and Max Lieblich.
\newblock Bridgeland-stable moduli spaces for {K}-trivial surfaces, 2007.
\newblock arXiv:0708.2247.

\bibitem[AD02]{Aspinwall-Douglas:stability}
Paul~S. Aspinwall and Michael~R. Douglas.
\newblock D-brane stability and monodromy.
\newblock {\em J. High Energy Phys.}, 5(5):no. 31, 35, 2002.
\newblock arXiv:hep-th/0110071.

\bibitem[AL01]{Aspinwall-Lawrence:DC-zero-brane}
Paul~S. Aspinwall and Albion Lawrence.
\newblock Derived categories and zero-brane stability.
\newblock {\em J. High Energy Phys.}, 8(8):Paper 4, 27, 2001.
\newblock arXiv:hep-th/0104147.

\bibitem[AP06]{Abramovich-Polishchuk:t-structures}
Dan Abramovich and Alexander Polishchuk.
\newblock Sheaves of {$t$}-structures and valuative criteria for stable
  complexes.
\newblock {\em J. Reine Angew. Math.}, 590:89--130, 2006.
\newblock arXiv:math/0309435.

\bibitem[Asp05]{Aspinwall:Dbranes-CY}
Paul~S. Aspinwall.
\newblock D-branes on {C}alabi-{Y}au manifolds.
\newblock In {\em Progress in string theory}, pages 1--152. World Sci. Publ.,
  Hackensack, NJ, 2005.
\newblock arXiv:hep-th/0403166.

\bibitem[Bay09]{large-volume}
Arend Bayer.
\newblock Polynomial {B}ridgeland stability conditions and the large volume
  limit.
\newblock {\em Geom. Topol.}, 13(4):2389--2425, 2009.
\newblock arXiv:0712.1083.

\bibitem[BBMT11]{BBMT:Fujita}
Arend Bayer, Aaron Bertram, Emanuele Macr{\`{\i}}, and Yukinobu Toda.
\newblock Bridgeland stability conditions on threefolds {II}: {A}n application
  to {F}ujita's conjecture, 2011.
\newblock arXiv:1106.3430.

\bibitem[BM02]{Bridgeland-Maciocia:K3fibrations}
Tom Bridgeland and Antony Maciocia.
\newblock Fourier-{M}ukai transforms for {$K3$} and elliptic fibrations.
\newblock {\em J. Algebraic Geom.}, 11(4):629--657, 2002.
\newblock arXiv:math/9908022.

\bibitem[BM11]{localP2}
Arend Bayer and Emanuele Macr{\`{\i}}.
\newblock The space of stability conditions on the local projective plane.
\newblock {\em Duke Math. J.}, 160(2):263--322, 2011.
\newblock arXiv:0912.0043.

\bibitem[Bog78]{Bogomolov:Ineq}
F.~A. Bogomolov.
\newblock Holomorphic tensors and vector bundles on projective manifolds.
\newblock {\em Izv. Akad. Nauk SSSR Ser. Mat.}, 42(6):1227--1287, 1439, 1978.

\bibitem[Bri06]{Bridgeland:stab-CY}
Tom Bridgeland.
\newblock Stability conditions on a non-compact {C}alabi-{Y}au threefold.
\newblock {\em Comm. Math. Phys.}, 266(3):715--733, 2006.
\newblock arXiv:math/0509048.

\bibitem[Bri07]{Bridgeland:Stab}
Tom Bridgeland.
\newblock Stability conditions on triangulated categories.
\newblock {\em Ann. of Math. (2)}, 166(2):317--345, 2007.
\newblock arXiv:math/0212237.

\bibitem[Bri08]{Bridgeland:K3}
Tom Bridgeland.
\newblock Stability conditions on {$K3$} surfaces.
\newblock {\em Duke Math. J.}, 141(2):241--291, 2008.
\newblock arXiv:math/0307164.

\bibitem[Dou02]{Douglas:stability}
Michael~R. Douglas.
\newblock Dirichlet branes, homological mirror symmetry, and stability.
\newblock In {\em Proceedings of the {I}nternational {C}ongress of
  {M}athematicians, {V}ol. {III} ({B}eijing, 2002)}, pages 395--408, Beijing,
  2002. Higher Ed. Press.
\newblock arXiv:math/0207021.

\bibitem[DRY06]{DRY:attractors}
Michael~R. Douglas, Rene Reinbacher, and Shing-Tung Yau.
\newblock Branes, bundles and attractors: Bogomolov and beyond, 2006.
\newblock arXiv:math/0604597.

\bibitem[Gie79]{Gieseker:Bog}
D.~Gieseker.
\newblock On a theorem of {B}ogomolov on {C}hern classes of stable bundles.
\newblock {\em Amer. J. Math.}, 101(1):77--85, 1979.

\bibitem[Har80]{Harris:space-curves}
Joe Harris.
\newblock The genus of space curves.
\newblock {\em Math. Ann.}, 249(3):191--204, 1980.

\bibitem[Har94]{Hartshorne:space-curves}
Robin Hartshorne.
\newblock The genus of space curves.
\newblock {\em Ann. Univ. Ferrara Sez. VII (N.S.)}, 40:207--223 (1996), 1994.

\bibitem[HL10]{HL:Moduli}
Daniel Huybrechts and Manfred Lehn.
\newblock {\em The geometry of moduli spaces of sheaves}.
\newblock Cambridge Mathematical Library. Cambridge University Press,
  Cambridge, second edition, 2010.

\bibitem[HMS08]{HMS:generic_K3s}
Daniel Huybrechts, Emanuele Macr{\`{\i}}, and Paolo Stellari.
\newblock Stability conditions for generic {$K3$} categories.
\newblock {\em Compos. Math.}, 144(1):134--162, 2008.
\newblock arXiv:math/0608430.

\bibitem[HRO96]{Happel-al:tilting}
Dieter Happel, Idun Reiten, and Smal{\o}Sverre O.
\newblock Tilting in abelian categories and quasitilted algebras.
\newblock {\em Mem. Amer. Math. Soc.}, 120(575):viii+ 88, 1996.

\bibitem[IUU10]{Ishii-Ueda-Uehara}
Akira Ishii, Kazushi Ueda, and Hokuto Uehara.
\newblock Stability conditions on {$A_n$}-singularities.
\newblock {\em J. Differential Geom.}, 84(1):87--126, 2010.
\newblock arXiv:math/0609551.

\bibitem[Kap88]{Kapranov:flagvarieties}
M.~M. Kapranov.
\newblock On the derived categories of coherent sheaves on some homogeneous
  spaces.
\newblock {\em Invent. Math.}, 92(3):479--508, 1988.

\bibitem[Kon95]{KontsICM94}
Maxim Kontsevich.
\newblock Homological algebra of mirror symmetry.
\newblock In {\em Proceedings of the {I}nternational {C}ongress of
  {M}athematicians, {V}ol.\ 1, 2 ({Z}\"urich, 1994)}, pages 120--139, Basel,
  1995. Birkh\"auser.
\newblock arXiv:alg-geom/9411018.

\bibitem[KS08]{Kontsevich-Soibelman:stability}
Maxim Kontsevich and Yan Soibelman.
\newblock Stability structures, motivic {D}onaldson-{T}homas invariants and
  cluster transformations, 2008.
\newblock arXiv:0811.2435.

\bibitem[Lan09a]{Langer:Survey}
Adrian Langer.
\newblock Moduli spaces of sheaves and principal {$G$}-bundles.
\newblock In {\em Algebraic geometry---{S}eattle 2005. {P}art 1}, volume~80 of
  {\em Proc. Sympos. Pure Math.}, pages 273--308. Amer. Math. Soc., Providence,
  RI, 2009.

\bibitem[Lan09b]{Langer:S-group}
Adrian Langer.
\newblock On the {S}-fundamental group scheme, 2009.
\newblock arXiv:0905.4600.

\bibitem[Mac04]{Macri:stability-examples}
Emanuele Macr{\`i}.
\newblock Some examples of moduli spaces of stability conditions on derived
  categories, 2004.
\newblock arXiv:math/0411613.

\bibitem[Mac07]{Macri:curves}
Emanuele Macr{\`{\i}}.
\newblock Stability conditions on curves.
\newblock {\em Math. Res. Lett.}, 14(4):657--672, 2007.
\newblock arXiv:0705.3794.

\bibitem[Mei07]{Sven:generic_tori}
Sven Meinhardt.
\newblock Stability conditions on generic complex tori, 2007.
\newblock arXiv:0708.3053.

\bibitem[Oka06]{Okada:P1}
So~Okada.
\newblock Stability manifold of {${\Bbb P}^1$}.
\newblock {\em J. Algebraic Geom.}, 15(3):487--505, 2006.
\newblock arXiv:math/0411220.

\bibitem[OS85]{Okonek-Spindler:spektrum-II}
Christian Okonek and Heinz Spindler.
\newblock Das {S}pektrum torsionsfreier {G}arben. {II}.
\newblock In {\em Seminar on deformations (\L\'od\'z/{W}arsaw, 1982/84)},
  volume 1165 of {\em Lecture Notes in Math.}, pages 211--234. Springer,
  Berlin, 1985.

\bibitem[Pol07]{Polishchuk:families-of-t-structures}
A.~Polishchuk.
\newblock Constant families of {$t$}-structures on derived categories of
  coherent sheaves.
\newblock {\em Mosc. Math. J.}, 7(1):109--134, 167, 2007.
\newblock arXiv:math/0606013.

\bibitem[Rei78]{Reid:Bog}
Miles Reid.
\newblock Bogomolov's theorem {$c_{1}^{2}\leq 4c_{2}$}.
\newblock In {\em Proceedings of the {I}nternational {S}ymposium on {A}lgebraic
  {G}eometry ({K}yoto {U}niv., {K}yoto, 1977)}, pages 623--642, Tokyo, 1978.
  Kinokuniya Book Store.

\bibitem[Sch80]{Schneider:Pn-rk3}
Michael Schneider.
\newblock Chernklassen semi-stabiler {V}ektorraumb\"undel vom {R}ang {$3$} auf
  dem komplex-projektiven {R}aum.
\newblock {\em J. Reine Angew. Math.}, 315:211--220, 1980.

\bibitem[Sim92]{Simpson:Higgs}
Carlos~T. Simpson.
\newblock Higgs bundles and local systems.
\newblock {\em Inst. Hautes \'Etudes Sci. Publ. Math.}, 75(75):5--95, 1992.

\bibitem[Tod08]{Toda:stab-crepant_res}
Yukinobu Toda.
\newblock Stability conditions and crepant small resolutions.
\newblock {\em Trans. Amer. Math. Soc.}, 360(11):6149--6178, 2008.
\newblock arXiv:math/0512648.

\bibitem[Tod09a]{Toda:limit-stable}
Yukinobu Toda.
\newblock Limit stable objects on {C}alabi-{Y}au 3-folds.
\newblock {\em Duke Math. J.}, 149(1):157--208, 2009.
\newblock arXiv:0803.2356.

\bibitem[Tod09b]{Toda:CY-fibrations}
Yukinobu Toda.
\newblock Stability conditions and {C}alabi-{Y}au fibrations.
\newblock {\em J. Algebraic Geom.}, 18(1):101--133, 2009.
\newblock arXiv:math/0608495.

\bibitem[Tod10]{Toda:PTDT}
Yukinobu Toda.
\newblock Curve counting theories via stable objects {I}. {DT}/{PT}
  correspondence.
\newblock {\em J. Amer. Math. Soc.}, 23(4):1119--1157, 2010.
\newblock arXiv:0902.4371.

\bibitem[Tod12]{Toda:ExtremalContractions}
Yukinobu Toda.
\newblock Stability conditions and extremal contractions, 2012.
\newblock arXiv:1204.0602.

\end{thebibliography}
\bibliographystyle{alphaspecial}     

\end{document}